\documentclass{article}

\usepackage{latexsym}
\usepackage{a4wide}
\usepackage{amscd}
\usepackage{graphics}
\usepackage{amsmath}
\usepackage{amssymb}
\usepackage{mathrsfs}
\input xy
\xyoption{all}

\newcommand{\N}{\mathbb{N}}                     
\newcommand{\Z}{\mathbb{Z}}                     
\newcommand{\R}{\mathbb{R}}                     
\newcommand{\C}{\mathbb{C}}                     
\newcommand{\set}[2]{\left\{{#1}\mid{#2}\right\}}       
\newcommand{\proof}{{\sl Proof.}\hspace{5pt}}   
\newcommand{\qed}{\hfill $\Box$ \bigskip}       
\newcommand{\dist}{\mathrm{dist\,}}             
\newcommand{\Ker}{\mathrm{Ker\,}}               
\newcommand{\coker}{\mathrm{coker\,}}           
\newcommand{\ind}{\mathrm{ind\,}}               
\newcommand{\codim}{\mathrm{codim}}           
\newcommand{\graf}{\mathrm{graph\,}}            
\newcommand{\ran}{\mathrm{ran\,}}		
\newcommand{\Det}{\mathrm{Det}}                 

\newtheorem{thm}{\sc Theorem}[section]      
\newtheorem{cor}[thm]{\sc Corollary}        
\newtheorem{lem}[thm]{\sc Lemma}            
\newtheorem{prop}[thm]{\sc Proposition}     
\newtheorem{defn}[thm]{\sc Definition}      
\newtheorem{rem}[thm]{\sc Remark}	    

\numberwithin{equation}{section}

\title{Infinite dimensional Grassmannians}
\author{ Alberto Abbondandolo \hspace{3cm} \hfill Pietro Majer 
\\ Dipartimento di Matematica \hspace{3cm} \hfill Dipartimento di Matematica
\\Universit\`a di Pisa \hfill Universit\`a di Pisa\\ 
Largo Bruno Pontecorvo, 5 \hfill Largo Bruno Pontecorvo, 5\\
56127 Pisa, Italy. \hfill 56127 Pisa, Italy.}

\begin{document}

\renewcommand{\theenumi}{\roman{enumi}}
\renewcommand{\labelenumi}{(\theenumi)}

\maketitle

\nocite{shu96}
\nocite{cpr93}
\nocite{zha94}

\section*{Introduction}

Infinite dimensional Grassmannians play a relevant role in many fields
of mathematics. For instance, they appear as classifying spaces in the
homotopy theory of classical groups (see \cite{bot59}), they provide a
natural setting to study completely integrable systems and loop groups (see
\cite{sat81,sw85,ps86,gue97,arb02}). Furthermore, they turn out to be
significant in the Morse theory of infinite dimensional 
manifolds (see \cite{cjs95}, \cite{ama05},\cite{jos02}).

This paper is devoted to a systematic study of the analytic and
homotopic properties of some infinite dimensional Grassmannians. 

After having established some useful facts about
the set of symmetric idempotent elements of a Banach *-algebra, we recall the 
properties of the Grassmannian $\mathrm{Gr}(H)$ of all closed subspaces
of a real or complex Hilbert space $H$, and then we focus our analysis 
on the following topics. 

\medskip

  We study the space of Fredholm pairs $\mathrm{\mathrm{Fp}}(H)$, i.e.\ the set of
  pairs of closed linear subspaces $(V,W)$ of $H$ having finite
  dimensional intersection and finite codimensional sum, and
  its open and closed subset $\mathrm{\mathrm{Fp}}^*(H)$ consisting of the Fredholm
  pairs $(V,W)$ with $\dim V = \dim W=\infty$. 
  These are analytic Banach manifolds, and we show that each component
  of $\mathrm{\mathrm{Fp}}^*(H)$ has the homotopy
  type of $BO$ - in the real case - or $BU$ - in the complex case, the
  classifying spaces of the infinite orthogonal group, respectively of
  the infinite unitary group (see \cite{as69} for similar results about
  the homotopy type of the space of skew-adjoint Fredholm operators). 
  
\medskip

  We study the Grassmannian $\mathrm{Gr_c}(V,H)$ of all compact perturbations
  of the closed linear subspace $V\subset H$, i.e.\ of all the closed
  linear subspaces $W\subset H$ such that the orthogonal projection
  onto $W$ differs from the orthogonal projection onto $V$ by a
  compact operator. This space is called restricted Grassmannian by
  some authors (see \cite{cjs95}, where this space is considered with
  a weaker topology, or \cite{ps86}, where the ideal of compact
  operators is replaced by the ideal of Hilbert-Schmidt operators). 
  It has the structure
  of an analytic Banach manifold, but when $\dim V = \codim V=\infty$,
  as a subset of $\mathrm{Gr}(H)$ it is just a topological
  submanifold. In the latter case, each component of $\mathrm{Gr_c}(V,H)$  has the homotopy type of $BO$ - if $H$ is real - or of
  $BU$ - if $H$ is complex. A similar result is proved in \cite{ps86}
  - the Hilbert-Schmidt case is treated in that reference, but the
  difference turns out to be irrelevant from the point of view of
  homotopy theory. The proof we give here is somehow more direct. 

\medskip

  We study the essential Grassmannian $\mathrm{Gr_e}(H)$, i.e. the quotient
  of $\mathrm{Gr}(H)$ by the equivalence relation $V\sim W$ iff $V$ is a compact
  perturbation of $W$. This is also an analytic Banach manifold,
  isometric to the space of symmetric idempotent elements in the
  Calkin algebra, and its homotopy type is easily determined. 
  Quotient spaces by
  finer equivalence relations, taking the Fredholm index into account,
  are also introduced and studied. 
  
\medskip

  We define the determinant bundle over $\mathrm{\mathrm{Fp}}(H)$ as an analytic
  line bundle whose fiber at $(V,W)\in \mathrm{\mathrm{Fp}}(H)$ is the top degree component in the exterior
  algebra of the finite dimensional space $(V\cap W) \oplus
  (H/(V+W))^*$. The space of Fredholm pairs is in some sense the
  larger space on which the determinant bundle can be defined, and from
  this it is possible to deduce the
  construction of other line bundles, such as the determinant bundle
  over the space of Fredholm operators introduced by Quillen
  \cite{qui85}, or the determinant bundle over the Grassmannian of
  compact perturbations of a given infinite dimensional and infinite
  codimensional space (see \cite{sw85,ps86}). We also show that the composition of Fredholm operators can be lifted analytically to the Quillen determinant bundles, and that such a lift has a useful associativity property. These facts are used to show how the sum of a finite dimensional space to one of the elements in a Fredholm pair can be lifted analytically to the determinant bundles, and to prove associativity for such a construction. 
 
In the case of a real Hilbert space $H$, the determinant bundle over the space of Fredholm pairs $\mathrm{Fp}^*(H)$ induces a double covering of $\mathrm{Fp}^*(H)$, called its the orientation bundle. This is an explicit presentation of the universal covering  of $\mathrm{Fp}^*(H)$. The associativity property of determinant bundles implies that the orientation bundle can be used quite effectively to determine orientations of finite dimensional linear spaces which are found as transverse intersections of infinite dimensional ones. These facts are used in \cite{ama05} to associate a chain complex with integer coefficients to functionals defined on infinite dimensional manifolds, having critical points with infinite Morse index.     

\medskip

Two appendices conclude the paper. In the first one, we prove the existence of a continuous global section for a continuous linear surjective map between Banach spaces, a result which is used several times in the paper. In the second one, we prove the above mentioned  associativity for the composition of Fredholm operators lifted to the determinant bundles. 

\medskip

This work was partially written while the second author was visiting the Forschungsinstitut f\"ur
Mathematik of the ETH in Z\"urich. He wishes to thank the Institute for the very warm hospitality.
We wish to express our gratitude to Andrea Maffei for an enlightening
discussion about determinant bundles. 

\section{A few facts about Banach algebras}

Let $\mathcal{A}$ be a real or complex Banach *-algebra, that is a
real or complex Banach algebra $\mathcal{A}$ endowed with an
involution map $*$. The involution $x\mapsto x^*$ satisfies the
following properties:
\[ 
x^{**}=x, \quad (xy)^* = y^* x^*, \quad \|x^*\|=x^*, \quad (x+y)^*=x^*
+ y^*, \quad (\lambda x)^* = \overline{\lambda} x^*,
\]
where $\lambda$ is a real respectively a complex number. Let
$\mathrm{Sym}(\mathcal{A}) = \set{x\in \mathcal{A}}{x^*=x}$ and
$\mathrm{Skew}(\mathcal{A})= \set{x\in \mathcal{A}}{x^*=-x}$ be 
the subspaces of symmetric and skew-symmetric
elements of $\mathcal{A}$, so that $\mathcal{A} = \mathrm{Sym}(\mathcal{A})
\oplus \mathrm{Skew}(\mathcal{A})$. 
Denote by $U(\mathcal{A})$ the group of unitary elements of
$\mathcal{A}$, 
\[
U(\mathcal{A}) := \set{u\in \mathcal{A}}{u^* u = uu^*=1},
\]
and by
$P^s(\mathcal{A})$ the set of symmetric idempotent elements of
$\mathcal{A}$,
\[
P^s(\mathcal{A}) := \set{p\in \mathrm{Sym}(\mathcal{A})}{p^2=p}.
\]
The aim of this section is to prove the following facts.  

\begin{prop}
\label{banalg1}
The set $P^s(\mathcal{A})$ of symmetric idempotent elements of the
Banach *-algebra $\mathcal{A}$ is an analytic Banach submanifold of $\mathcal{A}$.   
\end{prop}

\begin{prop}
\label{banalg2}
Consider the analytic group action
\[
U(\mathcal{A}) \times P^s(\mathcal{A}) \rightarrow P^s(\mathcal{A}),
\quad (u,p) \mapsto upu^{-1}.
\]
The orbit of every $p\in P^s(\mathcal{A})$,
$U(\mathcal{A})\cdot p := \set{upu^{-1}}{u\in U(\mathcal{A})}$,
is open and closed in $P^s(\mathcal{A})$, and the map
\[
U(\mathcal{A}) \rightarrow U(\mathcal{A})\cdot p, \quad u \mapsto
upu^{-1},
\]
is an analytic principal $I(p)$-bundle, where $I(p)$ is the isotropy
subgroup of $p$, 
\[
I(p):= \set{u\in U(\mathcal{A})}{upu^{-1}=p}.
\]
\end{prop}

\begin{prop}
\label{banalg3}
If $\Phi:\mathcal{A} \rightarrow \mathcal{B}$ is a surjective
homomorphism of Banach *-algebras, its restriction to the space of symmetric
idempotent elements
\[
\Phi|_{P^s(\mathcal{A})} : P^s(\mathcal{A}) \rightarrow
P^s(\mathcal{B}),
\]
is a $C^0$ fiber bundle (possibly with a non-constant fiber). 
It is an analytic fiber bundle if and only if the linear
subspace $\Ker \Phi\cap \mathrm{Skew}(\mathcal{A})$ has a direct summand in $\mathcal{A}$.
\end{prop}

\begin{rem} The above results hold also in a Banach algebra setting,
  by dropping the symmetry requirement, and by replacing the
  group of unitary elements by the group of invertible elements. 
\end{rem}

\paragraph{Square roots of 1.}
An element $p\in \mathcal{A}$ is symmetric and idempotent if and only
if $2p-1$ is a symmetric square root of 1. This fact allows us
to deduce the above propositions by analogous statements for the space
\[
Q^s(\mathcal{A}) := \set{q\in \mathrm{Sym}(\mathcal{A})}{q^2=1}.
\]
Denote by $\|\cdot\|$ the norm of $\mathcal{A}$, and by $B_r(a)$
the open ball of radius $r$ centered in $a\in \mathcal{A}$. The set
$Q^s(\mathcal{A})$ is a closed subset of $U(\mathcal{A})$, in particular
every $q\in Q^s(\mathcal{A})$ has norm 1. The elements
$1$ and $-1$ are isolated in $Q^s(\mathcal{A})$. More precisely
\begin{equation}
\label{zero}
Q^s(\mathcal{A})\cap B_2(1) =\{1\}, \quad Q^s(\mathcal{A}) \cap B_2(-1)
=\{-1\}.
\end{equation}
Indeed, $q^2=1$ implies $2(1-q)=(1-q)^2$, whence $2\|1-q\|\leq
\|1-q\|^2$, and hence either $\|1-q\|=0$ or $\|1-q\|\geq 2$, proving
the first identity. The second identity follows, because
$-Q^s(\mathcal{A})=Q^s(\mathcal{A})$. 

Notice that every $q\in Q^s(\mathcal{A})$
determines a decomposition $\mathrm{Sym}(\mathcal{A}) = A^s(q) \oplus
C^s(q)$ into the symmetric anti-commutant and the
symmetric commutant of $q$:
\[
A^s(q) = \set{a\in \mathrm{Sym}(\mathcal{A})}{aq=-qa}, \quad 
C^s(q) = \set{a\in \mathrm{Sym}(\mathcal{A})}{aq=qa}.
\]
Indeed, the corresponding projections are
\begin{eqnarray*}
& \mathrm{Sym}(\mathcal{A}) \rightarrow A^s(q) , \quad & a \mapsto \frac{1}{2}
(a-qaq), \\
& \mathrm{Sym}(\mathcal{A}) \rightarrow C^s(q) , \quad & a \mapsto \frac{1}{2}
(a+qaq).
\end{eqnarray*}

\begin{lem}
Let $q\in Q^s(\mathcal{A})$. Then the map
\[
\phi_q : B_{\sqrt{3}/2} (0) \cap A^s(q) \rightarrow C^s(q), \quad x
\mapsto (1-x^2)^{\frac{1}{2}} q,
\]
is analytic, and $Q^s(\mathcal{A})$ is locally the graph of $\phi_q$:
\[
Q^s(\mathcal{A}) \cap [(B_{\sqrt{3}/2}(0)\cap A^s(q)) \times (B_{1/2}(q)
\cap C^s(q))] = \graf \phi_q.
\]
\end{lem}

\proof
Let $x\in B_{\sqrt{3}/2}(0)\cap A^s(q)$. Since $\sqrt{3}/2<1$, a
square root of $1-x^2$ is well-defined:
\[
z = (1-x^2)^{\frac{1}{2}} = \sum_{k=0}^{\infty} (-1)^k \binom{1/2}{k} 
x^{2k} = 1 -
\sum_{k=1}^{\infty} \Bigl| \binom{1/2}{k}\Bigr| x^{2k}.
\]
So $\phi_q$ is analytic on $B_{\sqrt{3}/2}(0)$. Moreover, $z$ is invertible and
\[
z^{-1} = (1-x^2)^{-\frac{1}{2}} = \sum_{k=0}^{\infty} (-1)^k
\binom{-1/2}{k} x^{2k} =
\sum_{k=0}^{\infty} \Bigl| \binom{-1/2}{k} \Bigr| x^{2k}.     
\]
Since the coefficients of the above series are positive, we obtain the
inequalities
\begin{eqnarray}
\label{uno}
\|1-z\| \leq 1-(1-\|x\|^2)^{\frac{1}{2}} < \frac{1}{2}, \\
\label{due}
\|z^{-1}\|\leq (1-\|x\|^2)^{-\frac{1}{2}} <2, \\
\label{tre}
\|1-z^{-1}\|\leq (1-\|x\|^2)^{-\frac{1}{2}}-1 < 1.
\end{eqnarray}
From (\ref{uno}) we have
\[
\|q-\phi_q(x)\|=\|q-zq\|\leq \|1-z\| \, \|q\| < \frac{1}{2}.
\]
Since $x$ anti-commutes with $q$, $z$ commutes with $q$, and so does
$\phi_q(x)=zq$. Therefore, $\phi_q$ maps $B_{\sqrt{3}/2}(0)\cap
A^s(q)$ into $B_{1/2}(q)\cap C^s(q)$. 

Let $p\in Q^s(\mathcal{A})$, $p=x+y$ with $x\in B_{\sqrt{3}/2}(0)\cap
A^s(q)$, $y\in B_{1/2}(q)\cap C^s(q)$. We claim that $y=\phi_q(x)$,
that is $z^{-1}yq=1$, with $z$ as before. By (\ref{zero}), it is
enough to show that $(z^{-1}yq)^2=1$ and $\|z^{-1}yq-1\|<2$. Since
\[
x^2+y^2 \in C^s(q), \quad xy+yx\in A^s(q),
\]
and $x^2+y^2+xy+yx=p^2=1\in C^s(q)$, from $\mathrm{Sym}(\mathcal{A})=A^s(q)
\oplus C^s(q)$ we deduce that $x^2+y^2=1$ and $xy+yx=0$. Therefore,
$[y,x^2]=0$, whence $[y,z^{-1}]=0$. Using also $[q,z^{-1}]=0$, we
obtain that $(z^{-1}yq)^2=1$. Moreover,
\[
\|z^{-1} yq-1\| \leq \|z^{-1} yq - z^{-1}\| + \|z^{-1} -1\| \leq
\|z^{-1}\|\, \|yq-q^2\| + \|z^{-1} -1\| \leq \|z^{-1}\| \, \|y-q\| +
\|z^{-1}-1\|,
\]
so by (\ref{due}) and (\ref{tre}),
\[
\|z^{-1} yq-1\| \leq 2 \|y-q\|+ 1 <2,
\]
concluding the proof of the claim, by (\ref{zero}).

On the other hand, let $p=x+\phi_q(x)$, with $x\in
B_{\sqrt{3}/2}(0)\cap A^s(q)$. Then, since $xq=-qx$ and $zq=qz$, we
conclude that $p^2=1$.
\qed

The above lemma has the following consequence, which implies
Proposition \ref{banalg1}.

\begin{prop}
The set $Q^s(\mathcal{A})$ of symmetric square roots of $1$ in the
Banach *-algebra $\mathcal{A}$ is an analytic submanifold of
$\mathcal{A}$, and its tangent space at $q\in Q^s(\mathcal{A})$ is the
symmetric anti-commutant $A^s(q)$ of $q$.
\end{prop}

\paragraph{The action of $\mathbf{U(\mathcal{A})}$ on 
$\mathbf{Q^s(\mathcal{A})}$.}
The group $U(\mathcal{A})$ acts analytically on $Q^s(\mathcal{A})$,
by
\[
U(\mathcal{A}) \times Q^s(\mathcal{A}) \rightarrow Q^s(\mathcal{A}),
\quad (u,q) \mapsto uqu^{-1},
\]
and the isotropy subgroup of $q\in Q^s(\mathcal{A})$ is
\[
I(q) := \set{u\in U(\mathcal{A})}{uqu^{-1}=q} = U(\mathcal{A}) \cap
C(q),
\]
where $C(q)$ denotes the commutant subalgebra of the element $q$. The
bijection $P^s(\mathcal{A})\rightarrow Q^s(\mathcal{A})$, $p\mapsto
2p-1$, commutes with the respective actions of $U(\mathcal{A})$, so
Proposition \ref{banalg2} is implied by the following:

\begin{prop}
\label{lasta}
The orbit of every $q\in Q^s(\mathcal{A})$,
$U(\mathcal{A})\cdot q := \set{uqu^{-1}}{u\in U(\mathcal{A})}$,
is open and closed in $Q^s(\mathcal{A})$, and the map
\[
\psi_q: U(\mathcal{A}) \rightarrow U(\mathcal{A})\cdot q, \quad u \mapsto
uqu^{-1},
\]
is an analytic principal $(U(\mathcal{A})\cap C(q))$-bundle.
\end{prop}

\proof
By a well known simple argument (see for instance \cite[section 7.4]{ste51}), it is enough to
prove that the map $\psi_q$ has an analytic local section
near $q$.
Let $q\in Q^s(\mathcal{A})$. Note that $B_1(q)$ consists of invertible
elements: if $a\in B_1(q)$, then $\|aq-1\|=\|aq-q^2\|\leq \|a-q\|<1$,
thus $aq$ is invertible, and so is $a$.
If $p\in B_2(q)\cap Q^s(\mathcal{A})$, then $(p+q)/2\in B_1(q)$. So
$p+q$ is invertible, and the map
\[
s:B_2(q)\cap Q^s(\mathcal{A}) \rightarrow U(\mathcal{A}) , \quad p
\rightarrow (p+q)|p+q|^{-1} q
\]
is analytic. Since $q$ commutes with $(p+q)^2$, it commutes also with
$|p+q|^{-1}$, so
\[
ps(p)=(p+q)q|p+q|^{-1}q=(p+q)|p+q|^{-1} = s(p)q,
\]
implying that $s$ is an analytic local section near $q$ for the map
$\psi_q$.
\qed


\paragraph{Surjective homomorphisms.}
Since the bijection $P^s(\mathcal{A})\rightarrow
Q^s(\mathcal{A})$, $p\mapsto 2p-1$, commutes with Banach *-algebra
homomorphisms, statement (ii) of the following proposition implies 
Proposition \ref{banalg3}.

\begin{prop}
\label{last}
Let $\Phi:\mathcal{A} \rightarrow \mathcal{B}$ be a surjective
homomorphism of Banach *-algebras.
\begin{enumerate}
\item The group homomorphism
\[
\Phi_U := \Phi|_{U(\mathcal{A})} : U(\mathcal{A}) \rightarrow U(\mathcal{B})
\]
is a $C^0$ principal bundle (possibly with a non constant fiber).
\item The restriction
\[
\Phi_{Q^s}:=\Phi|_{Q^s(\mathcal{A})} : Q^s(\mathcal{A}) \rightarrow
Q^s(\mathcal{B}),
\]
is a $C^0$ fiber bundle (possibly with a non-constant fiber).
\item The maps $\Phi_U$ and $\Phi_{Q^s}$ are compatible with the group
actions, meaning that the diagram 
\begin{equation*}
\begin{CD}
U(\mathcal{A}) \times Q^s(\mathcal{A}) @>{\Phi_U \times \Phi_{Q^s}}>>
U(\mathcal{B}) \times Q^s(\mathcal{B}) \\ @VVV @VVV \\
Q^s(\mathcal{A}) @>{\Phi_{Q^s}}>> Q^s(\mathcal{B})
\end{CD}
\end{equation*}
commutes.
\end{enumerate}
The bundles in (i) and (ii) are analytic if and only if the linear
subspace $\Ker \Phi\cap \mathrm{Skew}(\mathcal{A})$ has a direct summand in
$\mathcal{A}$. 
\end{prop} 

\proof
The restriction of $\Phi$ to the subspaces of skew-symmetric elements,
\[
\Phi|_{\mathrm{Skew}(\mathcal{A})} : \mathrm{Skew}(\mathcal{A}) \rightarrow
\mathrm{Skew}(\mathcal{B}), 
\]
is surjective, because $\mathcal{A}$ and $\mathcal{B}$ are the
direct sum of the subspaces of their symmetric and skew-symmetric
elements, and because $\Phi$ is surjective. 
By Proposition \ref{propa} of Appendix A, the above map has a continuous
global section, which can be chosen to be linear if and only if the
linear subspace $\Ker \Phi \cap \mathrm{Skew}(\mathcal{A})$ has a direct
summand in $\mathrm{Skew}(\mathcal{A})$, or equivalently in
$\mathcal{A}$. The map $\Phi|_{\mathrm{Skew}(\mathcal{A})}$ is locally
conjugated at 0 to $\Phi_U$ 
by the exponential map. Hence the map $\Phi_U$ has a
continuous local section mapping $1_{\mathcal{B}}$ to
$1_{\mathcal{A}}$. As already mentioned, this is equivalent to (i).

Let $q\in Q^s(\mathcal{B})$. By Proposition \ref{lasta}, the map
\[
\psi_q : U(\mathcal{B}) \rightarrow U(\mathcal{B}) \cdot q, \quad u
\mapsto uqu^{-1},
\]
has an analytic local section mapping $q$ to $1_{\mathcal{B}}$, so by
composition we get a continuous local section for the map
$\psi_q \circ \Phi_U$: there exists a neighborhood $\mathcal{N}\subset
U(\mathcal{B}) \cdot q$ of $q$, and a map
\[
s: \mathcal{N} \rightarrow U(\mathcal{A}), \quad s(q) =
1_{\mathcal{A}},
\]
such that $\Phi(s(p)) q \Phi(s(p))^{-1} = p$ for every $p\in
\mathcal{N}$. Therefore, the continuous map
\[
\mathcal{N} \times \Phi_{Q^s}^{-1}(\{q\}) \rightarrow \Phi_{Q^s}^{-1}
(\mathcal{N}), \quad (p,x) \mapsto s(p) x s(p)^{-1}.
\]
is a local trivialization for $\Phi_{Q^s}$ near $q$, proving (ii).

Claim (iii) is readily checked. The last statement follows from the
fact that $\Phi_U$ has an analytic local section if and only if 
$\ker \Phi \cap \mathrm{Skew} (\mathcal{A})$ has a direct summand in $\mathcal{A}$.
\qed

\begin{rem}
\label{remo}
The maps $\Phi_U$ and $\Phi_{Q^s}$ need not be surjective. Indeed, if
$\Phi:\mathcal{A} \rightarrow \mathcal{B}$ is the quotient map from the 
Banach algebra of bounded linear operators on a
Hilbert space onto its Calkin algebra, then $\Phi_U$ is not onto: 
indeed, the image of
$\Phi_U$ is the set of unitary elements in the Calkin algebra having
index zero. 
As for the map $\Phi_{Q^s}$, notice that if $X$ is a compact Hausdorff
space and $Y$ is a closed subspace, the map $\Phi:C(X)\rightarrow
C(Y)$, $f\mapsto f|_Y$, is surjective, but $\Phi_{Q^s}$ is surjective
if and only if each connected component of $X$ contains at most one
connected component of $Y$ (by the Gel'fand duality, see
\cite[Theorem 4.29]{dou98}, this situation covers the general case of
commutative complex $C^*$-algebras).  
\end{rem}

\section{The Hilbert Grassmannian and the space of Fredholm pairs}
\label{secuno}
 
By $\mathrm{L}(E,F)$, respectively $\mathrm{L_c}(E,F)$, we
denote the space of continuous linear, respectively compact linear, 
maps from the Banach space $E$ to the Banach space $F$. 
If $F=E$ we also use the short forms $\mathrm{L}(E)$ and $\mathrm{L_c}(E)$. 
The group of invertible linear maps on $E$ is denoted by $\mathrm{GL}(E)$,
and $\mathrm{GL_c}(E)$ denotes the subgroup of the compact perturbations of
the identity. The norm of the operator
$T\in \mathrm{L}(E,F)$ is denoted by $\|T\|$. 

Let now $H$ be a real or complex infinite dimensional separable Hilbert
space. The orthogonal projection onto a closed subspace
$V\subset H$ is denoted by $P_V$, while $V^{\perp}$ denotes
the orthogonal complement of $V$ in $H$. The unitary group of $H$ is denoted by $U(H)=\set{U\in \mathrm{GL}(H)}{U^*U=UU^*=I}$ (in the real case
it would be more appropriate to call $U(H)$ the orthogonal group, and
to denote it by $O(H)$, but since we wish to deal with the complex and
the real case simultaneously, we do not make such a distinction). 

Let $\mathrm{Gr}(H)$ be the {\em Grassmannian of} $H$, i.e.\ the set of
closed linear subspaces of $H$. 
The assignment $V\mapsto P_V$ is an
inclusion of $\mathrm{Gr}(H)$ into $\mathrm{L}(H)$, onto the closed
subset of the orthogonal projectors of $H$. We can therefore define,
for any $V,W\in \mathrm{Gr}(H)$ the distance
\[
\dist(V,W) := \|P_{V} - P_{W}\|,
\]
which makes $\mathrm{Gr}(H)$ a complete metric space. We always consider
the topology on $\mathrm{Gr}(H)$ induced by this distance. The weak topologies
of $\mathrm{L}(H)$ induce other interesting topologies on the Grassmannian, whose
properties are examined in \cite{shu96}.

\paragraph{Analytic structure of $\mathbf{\mathrm{Gr}(H)}$.} 
For any $V\in \mathrm{Gr}(H)$, the unit open ball $B_1(V)\subset \mathrm{Gr}(H)$ is
contractible, as shown by the homeomorphism
\[
\Psi_V : B_1(V) \rightarrow \mathrm{L}(V,V^{\perp}), \quad 
W \mapsto P_{V^{\perp}} (P_V|_W)^{-1},
\]
the inverse mapping being
\[
\Psi_V^{-1} : \mathrm{L}(V,V^{\perp}) \rightarrow B_1(V), \quad
A \mapsto \graf A.
\]
Indeed, $\Psi_V$ is well-defined because $P_V|_W$ is the restriction
of the isomorphism
\[
I-(P_W-P_V)(P_W - P_{W^{\perp}}) = P_V P_W + P_{V^{\perp}}
P_{W^{\perp}} : W \oplus W^{\perp} \rightarrow V \oplus V^{\perp}.
\]
The collection $\{\Psi_V\}_{V\in \mathrm{Gr}(H)}$ is an analytic
atlas, with transition maps
\[
\Psi_W \circ \Psi_V^{-1} (A) = P_{W^{\perp}} (I+A) [P_W (I+A)]^{-1}.
\]
Therefore, $\mathrm{Gr}(H)$ is an analytic Banach manifold. Actually,
the identification $V\mapsto P_V$ allows to see $\mathrm{Gr}(H)$ as an
analytic submanifold of $\mathrm{L}(H)$, by Proposition
\ref{banalg1}.

\begin{rem} It easily seen that the analytic structure on $\mathrm{Gr}(H)$ does
  not depend on the choice of the Hilbert inner product on $H$. In
  other words, if $H=V \oplus V^{\prime}$ then the map
\[
\mathrm{L}(V,V^{\prime}) \rightarrow \mathrm{Gr}(H), \quad A \mapsto \graf A,
\]
is an analytic coordinate system.
\end{rem}

\paragraph{Transverse intersections.} We often use the fact that
the map
\[
\set{(V,W) \in \mathrm{Gr}(H) \times \mathrm{Gr}(H)}{V+W=H}, \quad (V,W) \mapsto V \cap
W,
\]
is analytic. To prove this fact, let $(V,W)$ be a pair of closed
linear subspaces such that $V+W=H$ and set $H_0 := V\cap W$, $H_1 :=
H_0^{\perp} \cap V$, $H_2 := H_0^{\perp} \cap W$, so that
\[
H = H_0 \oplus H_1 \oplus H_2, \quad V = H_0 \oplus H_1, \quad W=H_0
\oplus H_2.
\]
Consider the analytic coordinate system
\[
\mathrm{L}(H_0 \oplus H_1,H_2) \times \mathrm{L}(H_0\oplus H_2,H_1)
\rightarrow \mathrm{Gr}(H) \times \mathrm{Gr}(H), \quad (A,B) \mapsto (\graf A, \graf
B),
\]
mapping $(0,0)$ into $(V,W)$. If $\|A\|<1$ and $\|B\|<1$, then
\[
\graf A \cap \graf B = \graf T(A,B),
\]
where $T(A,B)\in \mathrm{L}(H_0,H_1 \oplus H_2)$ is the linear
mapping
\begin{eqnarray*}
T(A,B) = (I_{H_0 \oplus H_1} -BA)^{-1} (B+BA) + (I_{H_0 \oplus
  H_2} - AB)^{-1} (A+AB) \\ = A+B+AB+BA+ABA+BAB+ABAB+BABA+ \dots,
\end{eqnarray*}
which depends analytically on $(A,B)$, proving the claim.

We may also express $P_{V\cap W}$ in terms of $P_V$ and $P_W$.
First notice that if $V+W$ is closed then 
\begin{equation}
\label{mdu}
\|P_V P_W - P_{V\cap W}\|<1.
\end{equation}
Indeed, in the case $H=V\oplus W$, $W$ is the graph of some $L\in
\mathrm{L}(V^{\perp}, V)$, and (\ref{mdu}) reduces to
\[
\|P_V P_W\| = \|P_V|_W\| = \frac{\|L\|}{\sqrt{1+\|L\|^2}} < 1.
\]
The general case of a closed sum follows because
\[
P_V P_W - P_{V\cap W} = P_{V^{\prime}} P_{W^{\prime}},
\]
where $V^{\prime}:= (V\cap W)^{\perp} \cap V$, $W^{\prime}:= (V\cap
W)^{\perp} \cap W$ are in direct sum in the closed linear subspace
$(V\cap W)^{\perp} \cap (V+W)$.
 
From the fact that $X\subset Y$ implies $P_X P_Y = P_Y P_X = P_X$ it
follows that 
\[
(P_V P_W)^n - P_{V\cap W} = (P_V P_W-P_{V\cap W})^n.
\]
Together with (\ref{mdu}), this implies that when $V+W$ is closed,
\[
P_{V\cap W} = \lim_{n\rightarrow \infty} (P_V P_W)^n,
\]
in the operator norm. Notice also
that the above limit is uniform on the set of those pairs $(V,W)$
with $\|P_V P_W - P_{V\cap W}\|< \theta<1$. Since the function
$(V,W)\mapsto \|P_V P_W - P_{V\cap W}\|$ is continuous on the open set
of pairs $(V,W)$ such that $V+W=H$, we conclude that on this space the
analytic function $P_{V\cap W}$ is the locally uniform limit of
the sequence of monomials $(P_V P_W)^n$.

\paragraph{Homotopy of $\mathbf{Gr(H)}$.} 
The connected components of $\mathrm{Gr}(H)$ are the subsets
\[
\mathrm{Gr}_{n,k}(H):= \set{V\in \mathrm{Gr}(H)}{\dim V=n, \; \codim V=k}, \quad n,k\in
\N\cup \{\infty\}, \; n+k = \infty.
\]
The unitary group $U(H)$ 
acts analytically and transitively on each
component $\mathrm{Gr}_{n,k}(H)$, and if $V\in \mathrm{Gr}_{n,k}(H)$, the map
\[
U(H) \rightarrow \mathrm{Gr}_{n,k}(H), \quad U \mapsto UV,
\]
defines an analytic principal $U(V)\times U(V^{\perp})$-bundle, by
Proposition \ref{banalg2} (see also
\cite[Corollary 8.1]{luf67}). In particular,
the above map is a fibration. It is a well known result of Kuiper's
that $U(H)$ is contractible when $H$ is an
infinite dimensional Hilbert space (see \cite{kui65}), 
so the exact homotopy sequence
associated to the above fibration yields the following isomorphisms
\[
\pi_i(\mathrm{Gr}_{n,\infty}(H)) \cong \pi_i(\mathrm{Gr}_{\infty,n}(H)) = \begin{cases}
  0, & \mbox{if } i=0, \\ \pi_{i-1}(\mathrm{GL}(n)), & \mbox{if } i\geq 1, 
  \end{cases} \quad \pi_i(\mathrm{Gr}_{\infty,\infty}(H))=0, \;\; \forall i\in
  \N.
\]
Therefore, a theorem by Whitehead (see the Corollary after Theorem 15
in \cite{pal66}) implies that $\mathrm{Gr}_{\infty,\infty}(H)$
is contractible, while $\mathrm{Gr}_{n,\infty}(H)$ and $\mathrm{Gr}_{\infty,n}(H)$ have
the homotopy type of the classifying space of $\mathrm{GL}(n)$.

\paragraph{Stiefel spaces.}
Let $n,k\in \N\cup \infty$, $n+k=\infty$. The {\em Stiefel spaces}
$\mathrm{St}_{n,k}(H)$ are the sets
\begin{eqnarray*}
\mathrm{St}_{n,\infty}(H) & := & \set{T\in \mathrm{L}(\R^n,H)}{T \mbox{ is
    injective}}, \quad n\in \N,\\
\mathrm{St}_{\infty,n}(H) & := & \set{T\in \mathrm{L}(H)}{T \mbox{ is 
semi-Fredholm of
    index } -n \mbox{ and injective}}, \quad n\in \N\cup \{\infty\}, 
\end{eqnarray*}
endowed with the operator norm topology. The map
\[
\mathrm{St}_{n,k}(H) \rightarrow \mathrm{Gr}_{n,k}(H), \quad T \mapsto \ran T,
\]
defines a $\mathrm{GL}(n)$-principal bundle, in the case $n<\infty$, and a
$\mathrm{GL}(H)$-principal bundle, in the case $n=\infty$ (see \cite[Theorem 8.6]{luf67}). The Stiefel spaces $\mathrm{St}_{n,\infty}(H)$ are contractible,
for every $n\in \N\cup \{\infty\}$, while
\[
\pi_i(\mathrm{St}_{\infty,n}(H)) = \begin{cases} 0, & \mbox{if } i=0, \\
\pi_{i-1}(\mathrm{GL}(n)), & \mbox{if } i\geq 1, \end{cases}
\]
for $n\in \N$ (see \cite[Corollary 8.4]{luf67}).

\paragraph{Fredholm pairs.}
A pair $(V,W)$ of closed subspaces of $H$ is said a {\em Fredholm
pair} if $V\cap W$ is finite dimensional, and $V+W$ is 
finite codimensional, hence closed (see \cite[section IV \S 4]{kat80}). 
In this situation, the {\em Fredholm index} of
  $(V,W)$ is the integer number
\[
\ind (V,W) = \dim V\cap W - \codim (V+W).
\]
The pair $(V,W)$ is Fredholm if and only if the operator
$P_{W^{\perp}}|_V\in \mathrm{L}(V,W^{\perp})$ is Fredholm, in which
case $\ind (V,W) = \ind (P_{W^{\perp}}|_V)$.

The set of Fredholm pairs in $H$ is denoted by $\mathrm{\mathrm{Fp}}(H)$.
If $X,Y\in \mathrm{Gr}(H)$, with $\dim X<\infty$ and $\codim Y<\infty$, then
$(X,Y)$ is a Fredholm pair, of index
\[
\ind (X,Y) = \dim X - \codim Y.
\]
If $(V,W)$ is a Fredholm pair, and $V\in \mathrm{Gr}_{\infty,\infty}(H)$, then
also $W\in \mathrm{Gr}_{\infty,\infty}(H)$. Therefore we set
\[
\mathrm{\mathrm{Fp}}^*(H) := \set{ (V,W) \in \mathrm{\mathrm{Fp}}(H)}{V,W\in \mathrm{Gr}_{\infty,\infty}(H)}.
\]

\paragraph{Homotopy of $\mathbf{Fp(H)}$.} 
The space $\mathrm{Fp}(H)$ is
open in $\mathrm{Gr}(H)\times \mathrm{Gr}(H)$, and the index is a continuous function on 
$\mathrm{Fp}(H)$. It is easily checked that the connected components of $\mathrm{Fp}(H)$ 
are the subsets
\begin{eqnarray*}
\mathrm{Gr}_{n,\infty}(H)\times \mathrm{Gr}_{\infty,m}(H), \quad
\mathrm{Gr}_{\infty,n}(H)\times \mathrm{Gr}_{m,\infty}(H), \quad n,m\in \N, \\
\mathrm{Fp}^*_k (H) :=\set{(V,W)\in \mathrm{Fp}^*(H)}{\ind (V,W)=k}, \quad k\in \Z.
\end{eqnarray*}
The homotopy type of the first two families of spaces can be deduced
from the results seen before. As for the homotopy type of $\mathrm{Fp}^*_k(H)$,
we recall some definitions and basic facts.      

Let $GL$ be  be the infinite linear group, that is the inductive limit
$GL=\displaystyle{\lim_{\rightarrow}} \,\mathrm{GL}(n)$. 
In the complex case, $GL$ has the homotopy type of
$U=\displaystyle{\lim_{\rightarrow}} \,U(n)$, the infinite unitary group, while
in the real case it has the homotopy type of
$O=\displaystyle{\lim_{\rightarrow}} \, O(n)$, the infinite orthogonal
group. 
By the Bott
periodicity theorem \cite{bot59}, for $i\geq 0$,
\begin{equation}
\label{bot1} 
\pi_i(U) = \begin{cases} \Z & \mbox{if $i$ is odd}, \\ 0 & \mbox{if
    $i$ is even}, \end{cases} \quad\quad 
\pi_i(O) = \begin{cases} \Z_2 & \mbox{if } i\equiv 0,1 \;
  (\mathrm{mod}\; 8), \\ 0  & \mbox{if } i\equiv 2,4,5,6 \;
  (\mathrm{mod}\; 8), \\ \Z  & \mbox{if } i\equiv 3,7 \;
  (\mathrm{mod}\; 8). \end{cases}  
\end{equation}
By $BGL$ we denote the classifying space (see \cite[\S 7]{dol63})
of the topological group $GL$: in the complex case, $BGL$ has the
homotopy type of $BU$, the classifying space of $U$, while in the real
case it has the homotopy type of $BO$, the classifying space of
$O$. The spaces $BU$ and $BO$ are path connected, and (\ref{bot1})
implies that for $i\geq 1$,
\[
\pi_i(BU) = \begin{cases} \Z & \mbox{if } i \mbox{ is even}, \\ 
0  & \mbox{if } i \mbox{ is odd,} \end{cases} \quad\quad
\pi_i(BO) = \begin{cases} \Z & \mbox{if } i\equiv 0,4 \;
  (\mathrm{mod}\; 8), \\ \Z_2  & \mbox{if } i\equiv 1,2 \;
  (\mathrm{mod}\; 8), \\ 0  & \mbox{if } i\equiv 3,5,6,7 \;
  (\mathrm{mod}\; 8). \end{cases}
\] 

\begin{thm}
\label{hhoomm}
The components $\mathrm{Fp}_k^*(H)$ are pairwise homeomorphic, and they have
the homotopy type of $BU$ - if $H$ is a complex space - or $BO$ - if
$H$ is a real space.
\end{thm}

The proof of this result is given at the end of section \ref{essgrass}.

\section{The Grassmannian of compact perturbations}

We say that the subspace $W$ is a {\em compact perturbation} of
$V$ if its orthogonal projector $P_W$ is a compact perturbation of $P_V$. It is an equivalence relation. 

If $V$ is a compact perturbation of $W$, then $(V,W^{\perp})$ is a Fredholm
pair, and the operator $P_W|_V:V\rightarrow W$ is Fredholm. 
The {\em relative dimension of $V$ with respect to
 $W$} is the integer 
\begin{equation}
\label{numo}
\dim(V,W) := \ind (V,W^{\perp}) 
= \dim V\cap W^{\perp} - \dim V^{\perp} \cap W = \ind
(P_W|_V:V\rightarrow W). 
\end{equation}
See \cite{ama01} for a proof of the identities above. See also \cite[Remark 4.9]{bdf73} for an early appearance of the concepts of compact perturbation and of relative dimension for linear subspaces.
When $V$ and $W$ are finite dimensional (resp.\ finite codimensional),
we have $\dim (V,W)=\dim V-\dim W$ (resp.\ $\dim (V,W)=\codim W -
\codim V$). It is easy to see that the image by a linear isomorphism $T$ of a compact perturbation of some $V\in \mathrm{Gr}(H)$ is a compact perturbation of $TV$, and that the relative dimension is preserved (see \cite{ama01},
Proposition 2.4). Hence nor the notion of compact perturbation, neither the relative dimension depend on the choice of an inner product in $H$.

\begin{prop}
\label{ugo}
If $(W,Z)$ is a Fredholm pair of subspaces and $V$ is a compact
perturbation of $W$, then $(V,Z)$ is a Fredholm pair, with
\[
\ind (V,Z) = \ind (W,Z) + \dim (V,W).
\]
\end{prop}

\proof
From the identity 
\begin{equation}
\label{idecomp}
P_{Z^{\perp}}|_V = P_{Z^{\perp}}|_W \circ P_W|_V + P_{Z^{\perp}} (P_V -
P_W)|_V,
\end{equation}
we see that $P_{Z^{\perp}}|_V\in \mathrm{L}(V,Z^{\perp})$ is a
compact perturbation of the Fredholm operator
\[
V \stackrel{P_W|_V}{\longrightarrow} W
\stackrel{P_{Z^{\perp}}|_W}{\longrightarrow} Z^{\perp}.
\]
Therefore $P_{Z^{\perp}}|_V$ is Fredholm, and the conclusion follows
from (\ref{numo}) and from the additivity of the index by composition.
\qed

In particular, if each of the subspaces $V,W,Y$ is a compact perturbation of the other subspaces, there holds
\[
\dim(Y,V) = \dim (Y,W) + \dim (W,V).
\]
The following fact is proved in \cite{ama01}, Proposition 2.3.

\begin{prop}
\label{p1}
Let $H_1,H_2$ be Hilbert spaces, let $T,T^{\prime}\in
\mathrm{L}(H_1,H_2)$ be such that $T^{\prime}$ is a compact
perturbation of $T$, $\ran
T$ and $\ran T^{\prime}$ are closed. Then $\ker T^{\prime}$ is
a compact perturbation of $\ker T$, $\ran T^{\prime}$ is a compact
perturbation of $\ran T$, and
\[
\dim(\ran T^{\prime},\ran T) = - \dim (\ker T^{\prime},\ker T).
\]
\end{prop}

If $V\in \mathrm{Gr}(H)$, the {\em Grassmannian of compact perturbations of} $V$,
\[
\mathrm{Gr_c}(V,H) := \set{W\in \mathrm{Gr}(H)}{W \mbox{ is a compact
perturbation of }V}
\]
is a closed subspace of $\mathrm{Gr}(H)$. If $V$ has finite dimension
(respectively finite codimension), then
\[
\mathrm{Gr_c}(V,H) = \bigcup_{n\in \N} \mathrm{Gr}_{n,\infty}(H), \quad \Bigl(\mbox{
    resp. } = \bigcup_{n\in \N} \mathrm{Gr}_{\infty,n}(H) \Bigr).
\]
The interesting case arises when $V$ has both infinite dimension and
infinite codimension. In such a situation, $\mathrm{Gr_c}(V,H)$ is a
closed proper subset of $\mathrm{Gr}_{\infty,\infty}(H)$.
Moreover, the continuity of the Fredholm index implies that the function
\[
\dim(\cdot,\cdot) : \mathrm{Gr_c}(V,H) \times \mathrm{Gr_c}(V,H)
\rightarrow \Z
\]
is continuous. So the subsets
\[
\mathrm{Gr}_{\mathrm{c},n}(V,H) := \set{W\in \mathrm{Gr_c}(V,H)}{\dim (W,V)=n}, \quad n\in \Z.
\]
are open and closed in $\mathrm{Gr_c}(V,H)$. 
It is easily checked that the above sets are
connected, so they are the connected components of $\mathrm{Gr_c}(V,H)$.

\paragraph{Analytic structure of $\mathbf{Gr_c(V,H)}$.}
The restriction of the
map $\Psi_V$ defined in section \ref{secuno}
to the unit open ball $B_1(V)\cap \mathrm{Gr_c}(V,H)$ is a homeomorphism
onto $\mathrm{L_c}(V,V^{\perp})$. The collection of these
homeomorphisms defines an analytic atlas on $\mathrm{Gr_c}(V,H)$, which is
therefore an analytic Banach manifold modeled on
$\mathrm{L_c}(V,V^{\perp})$. Actually, Proposition \ref{banalg1}
implies that $\mathrm{Gr}(V,H)$ can be seen as an analytic Banach submanifold
of the Banach affine space $P_V + \mathrm{L_c}(H)$.
However, if both $V$ and
$V^{\perp}$ are infinite-dimensional, $\mathrm{Gr_c}(V,H)$ is not a $C^1$-submanifold
of $\mathrm{Gr}(H)$, because in this case the subspace
$\mathrm{L_c}(V,V^{\perp})$ does not admit a direct summand
in $\mathrm{L}(V,V^{\perp})$. Indeed, if $H$ is an
infinite dimensional separable Hilbert space, $\mathrm{L_c}(H)$ is
not complemented in $\mathrm{L}(H)$, a fact which can be
deduced\footnote{In fact, let us fix an orthonormal
  basis $\{e_n\}_{n\in \N}$ in $H$. Consider the canonical
  inclusion $j$ of $\ell^{\infty}$ into the space of diagonal
  operators with respect to this basis, and its left inverse $r:
  \mathrm{L}(H) \rightarrow \ell^{\infty}$ mapping $T\in
  \mathrm{L}(H)$ into the sequence $(\langle Te_n,e_n \rangle)_{n\in
    \N}\in \ell^{\infty}$. Then $j(c_0)\subset \mathrm{L_c}(H)$, and
  $r(\mathrm{L_c}(H)) = c_0$. If by contradiction $P:\mathrm{L}(H)
  \rightarrow \mathrm{L}(H)$ is a linear projection onto
  $\mathrm{L_c}(H)$, the above properties of $j$ and $r$ imply that
  $PjrP=jrP$. It follows that $rPj:\ell^{\infty} \rightarrow
  \ell^{\infty}$ is a projection onto $c_0$.}
from the analogous and more known fact that $c_0$ is not 
complemented\footnote{\label{fnote}Actually,
  every infinite dimensional complemented subspace of
  $\ell^{\infty}$ is isomorphic to $\ell^{\infty}$, see \cite[Theorem 2.a.7]{lt77}.} in $\ell^{\infty}$.  
The fact
that every surjective continuous linear map between Banach spaces
admits a (in general non differentiable) global section (see 
Appendix A), implies that $\mathrm{Gr_c}(V,H)$ is a $C^0$-submanifold of $\mathrm{Gr}(H)$. 

\paragraph{The action of $\mathbf{GL_c(H)}$ on $\mathbf{Gr_c(V,H)}$.}
The group $\mathrm{GL_c}(H)$, consisting of the linear automorphisms of $H$
which are compact perturbations of the identity, is an open subset of
the affine Banach space $I+\mathrm{L_c}(H)$. It inherits therefore
the structure of an analytic group, although it is just a
$C^0$-subgroup of $\mathrm{GL}(H)$ (again, because $\mathrm{L_c}(H)$ does not
have a direct summand in $\mathrm{L}(H)$). By Proposition \ref{p1},
there is an analytic group action
\[
\mathrm{GL_c}(H) \times \mathrm{Gr_c}(V,H) \rightarrow \mathrm{Gr_c}(V,H), \quad (L,W) \mapsto LW,
\]
which preserves the connected components of $\mathrm{Gr_c}(V,H)$, and such an
action is easily seen to be transitive on each component. 
Fix some $V\in \mathrm{Gr}(H)$. We claim that the map
\begin{equation}
\label{anbun}
\mathrm{GL_c}(H) \rightarrow \mathrm{Gr}_{\mathrm{c},0}(V,H), \quad L \mapsto LV,
\end{equation}
is an analytic principal $G$-bundle, with $G=\set{L\in
\mathrm{GL_c}(H)}{LV=V}$. Indeed, an analytic local section of (\ref{anbun})
near $V$ is 
\[
s:B_1(V)\cap \mathrm{Gr}_{\mathrm{c},0}(V,H) \rightarrow \mathrm{GL_c}(H), \quad W \mapsto P_W P_V +
P_{W^{\perp}} P_{V^{\perp}} = I - (P_V - P_W)(P_V -P_{V^{\perp}}),   
\]
and such a local section defines a local trivialization near any
$W_0\in \mathrm{Gr}_{\mathrm{c},0}(V,H)$ in the standard way: fixing $L_0\in \mathrm{GL_c}(H)$ such
that $L_0 V = W_0$, the map
\[
L_0 [B_1(V)\cap \mathrm{Gr}_{\mathrm{c},0}(V,H)] \times G \rightarrow \mathrm{GL_c}(H), \quad
(W,L) \mapsto L_0 s(L_0^{-1}W)L,
\]
is an analytic trivialization of (\ref{anbun}) near $W_0$.

\paragraph{Homotopy of $\mathbf{Gr_c(V,H)}$.}
Let $V\in \mathrm{Gr}_{\infty,\infty}(H)$. The components $\mathrm{Gr}_{\mathrm{c},n}(V,H)$
are pairwise analytically diffeomorphic. More generally,
$\mathrm{Gr}_{\mathrm{c},0}(V,H)$ is analytically diffeomorphic to $\mathrm{Gr}_{\mathrm{c},0}(W,H)$, for $V,W\in
\mathrm{Gr}_{\infty,\infty}(H)$, a diffeomorphism being induced by an $L\in
\mathrm{GL}(H)$ mapping $V$ onto $W$. 

In order to study the homotopy type of $\mathrm{Gr}_{\mathrm{c},0}(V,H)$, we introduce the
{\em Stiefel space of compact perturbations of} $V\hookrightarrow H$,
\[
\mathrm{St_c}(V,H) := \set{T\in \mathrm{L}(V,E)}{T \mbox{ is injective, }
  Tx=x+Kx \mbox{ with } K\in \mathrm{L_c}(V,H)},
\]
endowed with the norm topology. The set $\mathrm{St_c}(V,H)$ is open in the affine Banach
space $I + \mathrm{L_c}(V,H)$, so it is an analytic manifold. The map
\begin{equation}
\label{pb}
\mathrm{St_c}(V,H) \rightarrow \mathrm{Gr}_{\mathrm{c},0}(V,H), \quad T \mapsto \ran T,
\end{equation}
is an analytic principal $\mathrm{GL_c}(V)$-bundle. Indeed, if $W_0\in
\mathrm{Gr}_{\mathrm{c},0}(V,H)$ and $s_0:\mathcal{U}\rightarrow \mathrm{GL_c}(H)$ is an analytic  local
section of the fiber bundle (\ref{anbun}), with $\mathcal{U}$ a
neighborhood of $W_0$ in $\mathrm{Gr}_{\mathrm{c},0}(V,H)$, the map
\[
\mathcal{U} \times \mathrm{GL_c}(V) \rightarrow \mathrm{St_c}(V,H), \quad (W,A) \mapsto s_0(W)A,
\]
is an analytic local trivialization of (\ref{pb}) near $W_0$.  

Let $(V_n)$ be an increasing sequence of
subspaces of $V$, with $\dim V_n=n$ and such that their union is dense in
$V$. Then the union of the closed subspaces 
\[
\mathcal{K}_n = \set{K\in \mathrm{L_c}(V,H)}{K P_{V_n^{\perp}} =0}, 
\quad n\in \N,
\]
is dense in $\mathrm{L_c}(V,H)$. A well
known result by Palais (see \cite{pal66}, Corollary after Theorem 17)  
implies that $\mathrm{St_c}(V,H)$ is homotopically
equivalent to the inductive limit of the spaces $\mathrm{St_c}(V,H) \cap
(I+\mathcal{K}_n)$, equipped with the limit topology.   
Moreover, the set $\mathrm{St_c}(V,H) \cap (I+\mathcal{K}_n)$ is homeomorphic to
the Stiefel space $\mathrm{St}_{n,\infty}(H)$, introduced in section
\ref{secuno}.  
Since such a space is contractible, this proves that $\mathrm{St_c}(V,H)$ is also
contractible. In particular, the principal bundle (\ref{pb}) is the
universal principal bundle for the topological group $\mathrm{GL_c}(V)$, and
$\mathrm{Gr}_{\mathrm{c},0}(V,H)$ is the classifying space for such a group. The exact
homotopy sequence now implies that
\[
\pi_i(\mathrm{Gr}_{\mathrm{c},0}(V,H)) \cong \pi_{i-1} (\mathrm{GL_c}(V)), \quad\mbox{for } i\geq 1.
\]
Palais \cite{pal65} has shown that $\mathrm{GL_c}(V)$ has the homotopy type of 
the infinite dimensional general linear group $GL$, hence of $U$ in the
complex case, and of $O$ in the real case. 
Therefore, we have proved the following result.

\begin{thm}
\label{hgpc}
Let $V\in \mathrm{Gr}_{\infty,\infty}(H)$. Then the components $\mathrm{Gr}_{\mathrm{c},n}(V,H)$ are
pairwise analytically diffeomorphic, and
$\pi_i(\mathrm{Gr}_{\mathrm{c},n}(V,H))\cong \pi_{i-1}(GL)$ for $i\geq 1$, so that $\mathrm{Gr}_{\mathrm{c},n}(V,H)$
has the homotopy type of $BU$ - if $H$ is complex - or $BO$ - if $H$
is real. 
\end{thm}

\section{Essential Grassmannians}
\label{essgrass}

We have seen that the notion of compact perturbation produces  an equivalence
relation on $\mathrm{Gr}(H)$. The {\em essential Grassmannian of} $H$ is the
space of the equivalence classes for such a relation, endowed with the
quotient topology, and it is denoted by $\mathrm{Gr_e}(H)$. Denote by
\[
\Pi : \mathrm{Gr}(H) \rightarrow \mathrm{Gr_e}(H)
\]
the quotient projection. 
This space can be described in terms of the 
Calkin algebra of $H$, i.e.\ the quotient $C^*$-algebra
$\mathcal{C}(H):= \mathrm{L}(H)/\mathrm{L_c}(H)$. Denote by $\pi:
\mathrm{L}(H) \rightarrow \mathcal{C}(H)$ the quotient projection.   
The closed subspace $V$ is a compact perturbation of $W$ if and only if
$\pi(P_V)=\pi(P_W)$, so $\pi$ maps the set of
orthogonal projectors into the set of symmetric idempotent elements of
the Calkin algebra. The next result shows that every symmetric idempotent
element of the Calkin algebra is the image of some orthogonal
projector (compare this result with Remark \ref{remo}). 

\begin{prop}
The restriction of the map $\pi$ to the set of the (orthogonal)
projectors of $H$ is onto the set of (symmetric) idempotent elements
of $\mathcal{C}(H)$.
\end{prop}

\proof In any Banach *-algebra, an element $x$ is (symmetric)
idempotent if and only if $2x-1$ is a (symmetric) square root of the
identity. It is therefore equivalent to show that the quotient
projector $\pi$ maps the (symmetric) square roots of $I$
in $\mathrm{L}(H)$ onto the set of (symmetric) square roots of $1$
in $\mathcal{C}(H)$. 

Let $Q\in \mathrm{L}(H)$ and $K\in \mathrm{L_c}(H)$ be such that
$Q^2=I-K$. We have to find $J\in \mathrm{L_c}(H)$ such that
$(Q-J)^2=I$. The operator $J$ will be self-adjoint if $Q$ and $K$ are.

The spectral mapping theorem implies that
\[
\sigma(Q)^2 = \sigma(Q^2) = \sigma(I-K),
\]
so the set $\sigma(Q)\setminus \{-1,1\}$ consists of isolated
eigenvalues of finite multiplicity. 

Let $U=\set{z\in \C}{|z^2-1|<1}$,
and let $H=H_0\oplus H_1$, $Q=Q_0 \oplus Q_1$ be the spectral
decomposition of the operator $Q$ corresponding to the decomposition
of $\sigma(Q)$ into the closed sets $\sigma(Q)\setminus U$ and
$\sigma(Q) \cap U$. Denote by $P_0,P_1$ the corresponding
projectors. Then $H_0$ is finite dimensional, and $Q_1$ is invertible,
with $\sigma(Q_1)\subset U$. Therefore, the spectrum of the operator
$K_1:=K|_{H_1}=I-Q_1^2\in \mathrm{L_c}(H_1)$ is contained in the unit
ball $\set{z\in \C}{|z|<1}$. The analytic function
\[
f(z) := 1-z-(1-z)^{\frac{1}{2}} = - \frac{z}{2} + \sum_{n=2}^{\infty}
\left|\binom{1/2}{n}\right| z^n, \quad |z|<1,
\]
verifies the identity
\begin{equation}
\label{ide}
\frac{f(z)^2}{1-z} - 2f(z) - z =0.
\end{equation}
Set $J_0:= Q_0 - P_0\in \mathrm{L}(H_0)$, 
$J_1:= Q_1^{-1} f(K_1)\in \mathrm{L_c}(H_1)$. Then $J:=J_0 \oplus
J_1$ is compact, and by (\ref{ide}) it satisfies $(Q-J)^2=I$.
\qed

The above result implies that $\pi$ induces a homeomorphism between the
essential Grassmannian $\mathrm{Gr_e}(H)$ and the space of idempotent
symmetric elements of $\mathcal{C}(H)$. Therefore $\mathrm{Gr_e}(H)$
inherits the structure of a complete metric space. Furthermore, 
Proposition \ref{banalg1} shows that $\mathrm{Gr_e}(H)$ can be given the
structure of an analytic Banach submanifold of $\mathcal{C}(H)$. 

\paragraph{Homotopy of $\mathbf{Gr_e(H)}$.}
The finite dimensional and the finite codimensional spaces represent
two isolated points in the essential Grassmannian. What remains is a
connected component, which is the image of $\mathrm{Gr}_{\infty,\infty}(H)$ by
$\Pi$, and which is denoted by $\mathrm{Gr_e^*}(H)$. By Proposition
\ref{banalg3}, the map
\[
\Pi: \mathrm{Gr}_{\infty,\infty}(H) \rightarrow \mathrm{Gr_e^*}(H)
\]
is a $C^0$ fiber bundle, with typical fiber $\mathrm{Gr_c}(V,H)$, for some
$V\in \mathrm{Gr}_{\infty,\infty}(H)$. In particular, the above map is a
fibration, and since its total space 
is contractible, the exact homotopy sequence implies that
\[
\pi_i(\mathrm{Gr_e^*}(H)) \cong \pi_{i-1}(\mathrm{Gr_c}(V,H)).
\]
Therefore, Theorem \ref{hgpc} has the following consequence.
\begin{cor}
The space $\mathrm{Gr_e^*}(H)$ is path connected, its fundamental group is
infinite cyclic, and
$\pi_i(\mathrm{Gr_e^*}(H))\cong \pi_{i-2}(GL)$ for $i\geq 2$.
\end{cor}

\paragraph{The $\mathbf{(m)}$-essential Grassmannian.}
In order to represent the coverings of $\mathrm{Gr_e^*}(H)$, we consider
the quotient of $\mathrm{Gr}(H)$ by stronger
equivalence relations, which take the relative dimension into account: 
if $m\in \N$, we define the $(m)$-essential Grassmannian $\mathrm{Gr}_{(m)}(H)$
to be the quotient of $\mathrm{Gr}(H)$ by the equivalence relation 
\[
\set{(V,W)\in \mathrm{Gr}(H)\times \mathrm{Gr}(H)}{V \mbox{ is a compact perturbation of
  } W, \mbox{ and } \dim (V,W)\in m\Z}.
\]
Denote by
\[
\Pi_m : \mathrm{Gr}_{\infty,\infty}(H) \rightarrow \mathrm{Gr}_{(m)}(H)
\]
the quotient map. With this terminology, the $(1)$-essential
Grassmannian is just the
essential Grassmannian, and the $(0)$-essential Grassmannian
distinguishes commensurable spaces according to their relative dimension. 
We can endow $\mathrm{Gr}_{(m)}(H)$ with the quotient topology induced
by $\mathrm{Gr}(H)$. 

Every set $\mathrm{Gr}_{n,\infty}(H)$ or $\mathrm{Gr}_{\infty,n}(H)$,
$n\in \N$, represents an isolated point in $\mathrm{Gr}_{(0)}(H)$, which has
thus infinitely many isolated points.
If $m\geq 1$, the sets
\[
\bigcup_{n\in m\Z + k} \mathrm{Gr}_{n,\infty}(H) \quad \mbox{and} \quad 
\bigcup_{n\in m\Z + k} \mathrm{Gr}_{\infty,n}(H), \quad k=0,1,\dots,m-1,
\]
represent distinct isolated points in $\mathrm{Gr}_{(m)}(H)$, which has thus
$2m$ isolated points.  

The remaining part of $\mathrm{Gr}_{(m)}(H)$ consists of the 
quotient of $\mathrm{Gr}_{\infty,\infty}(H)$, denoted by
$\mathrm{Gr}_{(m)}^*(H)$. Since $\mathrm{Gr}_{\infty,\infty}(H)$ is connected, so is
$\mathrm{Gr}_{(m)}^*(H)$, for any $m\in \N$.

\paragraph{The action of $\mathbf{U(H)}$ on $\mathbf{Gr_{(m)}^*(H)}$.}
The action of $U(H)$ on $\mathrm{Gr}_{\infty,\infty}(H)$ induces a transitive
action of $U(H)$ on the quotient $\mathrm{Gr}_{(m)}^*(H)$. This action presents
$\mathrm{Gr}_{(m)}^*(H)$ as a quotient of $U(H)$: fixing $V\in
\mathrm{Gr}_{\infty,\infty}(H)$, the map
\[
h_m : U(H) \rightarrow \mathrm{Gr}_{(m)}^*(H), \quad U \mapsto \Pi_m(UV),
\]
is indeed open and onto. Let us determine the isotropy subgroup
$h_m^{-1}([V])$. By the identity
\[
P_{UV} - P_V = U P_V U^* - P_V = [U,P_V] U^*, \quad \forall U\in U(H),
\]
the subspace $UV$ is a compact perturbation of $V$ if and only if $U$
belongs to the closed subgroup
\[
\mathcal{H} := \set{U\in U(H)}{[U,P_V]\in \mathrm{L_c}(H)}.
\]
If $U\in \mathcal{H}$, since
\[
P_V U P_V + P_{V^{\perp}} U P_{V^{\perp}} = U - [[U,P_V],P_V] \in U +
\mathrm{L_c}(H),
\]
the operator $P_V U P_V + P_{V^{\perp}} U P_{V^{\perp}}$ is Fredholm
of index 0. So $P_V U P_V$ and $P_{V^{\perp}} U P_{V^{\perp}}$ are
Fredholm operators on $V$ and $V^{\perp}$, respectively, with $\ind
(P_V U P_V) = -\ind (P_{V^{\perp}} U P_{V^{\perp}})$. By (\ref{numo}),
and by the additivity of the Fredholm index with respect to composition,
\begin{eqnarray*}
\dim(V,UV) = \ind (P_{UV}|_V:V\rightarrow UV) = \ind (U P_V
U^*|_V:V\rightarrow  UV) \\
= \ind (P_V U^*|_V: V\rightarrow V) + \ind (U|_V : V
\stackrel{\sim}{\rightarrow} UV) = \ind (P_V U^* P_V) = - \ind (P_V U
P_V).
\end{eqnarray*}
Therefore,
\[
\dim (UV,V) = \ind (P_V U P_V), \quad \forall U\in \mathcal{H},
\]
which shows that $h_m^{-1}([V])$ is the open and closed subgroup of
$\mathcal{H}$ ,
\[
\mathcal{H}_m := \set{U\in \mathcal{H}}{\ind (P_V U P_V) \in m\Z}.
\]
Therefore the space $\mathrm{Gr}_{(m)}^*(H)$ is homeomorphic to the space of
left cosets $U(H)/\mathcal{H}_m$. Now assume that $n\in \N^*$ divides
$m\in \N$. Since $\mathcal{H}_m$ is an open subgroup of
$\mathcal{H}_n$, the quotient map
\[
\Pi^m_n : \mathrm{Gr}_{(m)}^* (H) \rightarrow \mathrm{Gr}_{(n)}^*(H)
\]
is a covering, with fiber $n\Z/m\Z$. We claim that the map $h_m$ is a
$C^0$ principal $\mathcal{H}_m$-bundle. Since $h_1=\Pi^m_1 \circ h_m$,
it is enough to prove that $h_1$ is a
$C^0$ principal $\mathcal{H}_1$-bundle. In this case 
\[
\mathrm{Gr}_{(1)}^*(H)= \mathrm{Gr}_{e}^*(H) = P^s(\mathcal{C}) \setminus\{1,0\} =
\set{u\,\pi(P_V)\,u^{-1}}{u\in U(\mathcal{C}(H))},
\]
so by Proposition \ref{banalg2} the map
\[
U(\mathcal{C}(H)) \rightarrow \mathrm{Gr}_{(1)}^*(H), \quad u \mapsto u\,
\pi(P_V) \, u^{-1},
\]
has a continuous local section near $\Pi_1(V)$, mapping $\Pi_1(V)$ to $1$. 
By Proposition \ref{last}
(i), the quotient map $U(H) \rightarrow U(\mathcal{C}(H))$ has a
continuous local section near $1$. By composition we get a
continuous local section for $h_1$, which is therefore a $C^0$
principal $\mathcal{H}_1$-bundle by the standard trivialization. 

\paragraph{The homotopy of $\mathbf{Gr_{(m)}^*(H)}$.}
We claim that the quotient map
\[
\Pi_m|_{\mathrm{Gr}_{\infty,\infty}(H)} :\mathrm{Gr}_{\infty,\infty}(H) \rightarrow 
\mathrm{Gr}_{(m)}^*(H) 
\]
is a $C^0$ fiber bundle, with typical fiber
\[
\set{W\in \mathrm{Gr_c}(V,H)}{\dim(W,V)\in m\Z}.
\] 
Indeed, if $s:\mathcal{U} \subset 
\mathrm{Gr}_{(m)}^*(H)\rightarrow U(H)$ is a continuous local section for the
map $h_m$, the map
\[
\mathcal{U} \times \set{W\in \mathrm{Gr_c}(V,H)}{\dim(W,V)\in m\Z} \rightarrow
\mathrm{Gr}_{\infty,\infty}(H),\quad (\alpha,W) \mapsto s(\alpha)W,
\] 
is a local trivialization for the map $\Pi_m|_{\mathrm{Gr}_{\infty,\infty}(H)}$.

Since $\mathrm{Gr}_{\infty,\infty}(H)$ is contractible, by the exact homotopy
sequence of a fibration we immediately get the following result.

\begin{cor}
Let $m\in \N$. The space $\mathrm{Gr}_{(m)}^*(H)$ is
connected, its fundamental group is trivial for $m=0$, and
infinite cyclic for $m\geq 1$. Moreover, for every $m\in \N$,
\[
\pi_i(\mathrm{Gr}_{(m)}^*(H)) \cong \pi_i(\mathrm{Gr_e}(H)) \cong
\pi_{i-2}(GL), \quad \mbox{for }i\geq 2.
\]
\end{cor}

In particular, the quotient map
\[
\Pi_1^0: \mathrm{Gr}_{(0)}^* (H) \rightarrow \mathrm{Gr}_{(1)}^*(H) = \mathrm{Gr_e^*}(H)
\]
is the universal covering of $\mathrm{Gr_e^*}(H)$, and the
covering corresponding to the subgroup $m\Z\subset
\Z=\pi_1(\mathrm{Gr_e^*}(H))$, $m\in \N$, is the quotient map
\[
\Pi_1^m Gr^*_{(m)} (H) \rightarrow \mathrm{Gr}_{(1)}^*(H) = \mathrm{Gr_e^*}(H).
\]   
These are all the coverings of the essential Grassmannian.

We end this section by determining the homotopy groups of the
space of Fredholm pairs, thus proving Theorem \ref{hhoomm}.

\paragraph{Proof of Theorem \ref{hhoomm}.}
For $k\in \Z$, the map
\[
Q_k:\mathrm{Fp}_k^*(H) \rightarrow \mathrm{Gr}_{\infty,\infty}(H), \quad (V,W) \mapsto W,
\]
defines a fiber bundle. Indeed, if $W_0\in \mathrm{Gr}_{\infty,\infty}(H)$, the
map 
\[
B_1(W_0) \times Q_k^{-1}(\{W_0\}) \rightarrow \mathrm{Fp}_k^*(H), \quad (W,V)
\mapsto ((P_W P_{W_0} + P_{W^{\perp}} P_{W_0^{\perp}}) V,W),
\]
defines a local trivialization near $W_0$. Since
$\mathrm{Gr}_{\infty,\infty}(H)$ is contractible, $\mathrm{Fp}^*_k(H)$ is homeomorphic
to the product $\mathrm{Gr}_{\infty,\infty}(H)\times Q_k^{-1}(\{W_0\})$, for
any $W_0\in \mathrm{Gr}_{\infty,\infty}(H)$. Let $W_1$ be a compact
perturbation of $W_0$ with $\dim(W_0,W_1)=k$. Then by Proposition
\ref{ugo}, $Q_k^{-1}(\{W_0\})=Q_0^{-1}(\{W_1\})$, so $\mathrm{Fp}^*_k(H)$ is
homeomorphic to $\mathrm{Fp}^*_0(H)$, and it has the homotopy type of the fiber
of $Q_0$. 

There remains to determine the homotopy type of $Q_0^{-1}(\{W_0\})$,
for some fixed $W_0\in \mathrm{Gr}_{\infty,\infty}(H)$.
Consider the fiber bundle 
\[
\Pi: \mathrm{Gr}_{\infty,\infty}(H) \rightarrow \mathrm{Gr_e^*}(H),
\]
and the open subset of $\mathrm{Gr_e^*}(H)$
\[
\Omega := \set{\Pi(\graf T)}{T\in \mathrm{L}(W_0^{\perp},W_0)}.
\]
The space $\Omega$ is contractible: indeed the homotopy 
\[
\Omega \times [0,1] \rightarrow \Omega, \quad (\Pi(\graf T),\lambda)
\mapsto \Pi(\graf \lambda T)
\]
is well defined, because the graphs of $T,S\in
\mathrm{L}(W_0^{\perp},W_0)$ are commensurable if and only if $T-S$
is compact, by Proposition \ref{p1}. 
Moreover, using Proposition \ref{ugo} it is easy to see that
\[
\Pi^{-1}(\Omega) = \set{V\in \mathrm{Gr}_{\infty,\infty}(H)}{(V,W_0)\in \mathrm{Fp}(H)},
\]
and that $Q^{-1}(\{W_0\})$ is a connected component of
$\Pi^{-1}(\Omega)$, that is
\[
Q^{-1}(\{W_0\} = \set{V\in \mathrm{Gr}_{\mathrm{c},0}(\graf T,H)}{T\in
  \mathrm{L}(W_0^{\perp},W_0)}. 
\]
Therefore, the restriction of the map $\Pi$,
\[
\Pi_0 : Q^{-1}(\{W_0\}) \rightarrow \Omega, \quad V \mapsto \Pi(V),
\]
defines a fiber bundle, with typical fiber
\[
\Pi_0^{-1} (\Pi(W_0^{\perp})) = \mathrm{Gr}_{\mathrm{c},0}(W_0^{\perp},H).
\]
Since the base space of $\Pi_0$ is contractible, $Q_0^{-1}(\{W_0\})$
has the homotopy type of $\mathrm{Gr}_{\mathrm{c},0}(W_0^{\perp},H)$, and the conclusion
follows from Theorem \ref{hgpc}.
\qed 

\section{The functor $\mathbf{Det}$}
\label{Dets}

We denote by $\Det$ the functor which associates to any finite
dimensional real vector space the one-dimensional real vector space 
$\Lambda_{\max}(X)$, the component of top degree of the exterior
algebra of $X$, and to any linear map $T:X\rightarrow Y$ the linear map
\[
\Det(T) = \left\{ \begin{array}{ll} \Lambda_{\max}(T) & \mbox{if }
    \dim X = \dim Y, \\ 0 & \mbox{otherwise.} \end{array} \right.
\]
Therefore, $\Det(T)\neq 0$ if and only if $T$ is an isomorphism.

As it is well known, an exact sequence $T$ of finite dimensional vector spaces 
\[
0 \stackrel{T_0}{\longrightarrow} X_1 \stackrel{T_1}{\longrightarrow}
X_2 \stackrel{T_2}{\longrightarrow} \dots
\stackrel{T_{n-1}}{\longrightarrow} X_n
\stackrel{T_n}{\longrightarrow} 0
\]
induces an isomorphism
\[
\phi_T : \bigotimes_{i \;\mathrm{\scriptsize{ odd}}} \Det(X_i) \cong
\bigotimes_{i \;\mathrm{\scriptsize{ even}}} \Det(X_i).
\]
Actually, there are different conventions in the literature concerning the choice of the sign of such an isomorphism. As shown in the following sections, the choice of proper signs is important, so let us exhibit an explicit construction for the standard isomorphism $\phi_T$, and let us fix notations to denote the opposite isomorphism.

Consider an exact sequence $S$
\[
0 \stackrel{S_1}{\longleftarrow} X_1 \stackrel{S_2}{\longleftarrow}
X_2 \stackrel{S_3}{\longleftarrow} \dots
\stackrel{S_n}{\longleftarrow} X_n
\stackrel{S_{n+1}}{\longleftarrow} 0
\]
such that for every $i$
\begin{equation}
\label{inve}
T_{i-1} S_i + S_{i+1} T_i = I_{X_i}.
\end{equation}
Such an exact sequence is uniquely determined by the choice of an
algebraic linear complement $V_i$ for each subspace $\ker T_i$ in
$X_i$. In fact, $T_i$ restricts to an isomorphism from $V_i$ to $\ker
T_{i+1}$, and we can define $S_{i+1}:X_{i+1} \rightarrow X_i$ to
coincide with the inverse of such an isomorphism on $\ker T_{i+1}$ and
to be zero on $V_{i+1}$.

The linear maps $P_i := T_{i-1} S_i$ and $Q_i := S_{i+1} T_i$ are the
projectors associated to the splitting $X_i = \ker T_i \oplus \ker
S_i$. 

The linear map
\[
\Phi_T^S : \bigoplus_{i \;\mathrm{ \scriptsize{odd}}} X_i \longrightarrow
\bigoplus_{i \;\mathrm{ \scriptsize{even}}} X_i, \quad \Phi_T^S
=\bigoplus_{i \;\mathrm{ \scriptsize{odd}}} (T_i + S_i)
\]
is an isomorphism, its inverse being
\[
{\Phi_T^S}^{-1} =  \bigoplus_{i \;\mathrm{\scriptsize{ even}}} (T_i + S_i).
\]
If $S'$ is another exact sequence satisfying (\ref{inve}), then
\[
{\Phi_T^{S'}}^{-1} \circ \Phi_T^S = \bigoplus_{i \;\mathrm{\scriptsize{ odd}}}
(S'_{i+1} T_i + T_{i-1} S_i) = I +  \bigoplus_{i \;\mathrm{\scriptsize{
      odd}}} T_{i-1} (S_i - S_i'),
\]
and the last term is nilpotent,
\[
\left( \bigoplus_{i \;\mathrm{\scriptsize{odd}}} T_{i-1} (S_i - S_i')
\right)^2 = 0.
\]
It follows that
\[
\det {\Phi_T^{S'}}^{-1} \circ \Phi_T^S = 1,
\]
which implies
\[
\Det(\Phi_T^S) = \Det(\Phi_T^{S'}).
\]
Recalling that the exterior algebra of a direct sum is naturally
identified with the tensor product of the corresponding exterior
algebras, we can associate to the exact sequence $T$ the isomorphism
\[
\phi_T : \bigotimes_{i \;\mathrm{\scriptsize{ odd}}} \Det(X_i) \cong
\bigotimes_{i \;\mathrm{\scriptsize{ even}}} \Det(X_i), \quad
\phi_T := \Det(\Phi_T^S),
\]
where $S$ is any sequence satisfying (\ref{inve}). Such an isomorphism
is natural, in the sense that if we have an isomorphism of exact
sequences
\begin{equation}
\label{isoexa}
\begin{CD} 0 @>>> X_1 @>T_1>> X_2 @>T_2>> X_3 @>>> \dots @>>> X_n @>>>
  0 \\ & & @VR_1VV @VVR_2V @VVR_3V & & @VVV & \\ 0 @>>> X_1' @>T_1'>> X_2'
  @>T_2'>> X_3' @>>> \dots @>>> X_n' @>>> 0
\end{CD}
\end{equation}
there holds
\begin{equation}
\label{natural}
\left( \bigotimes_{i \;\mbox{\scriptsize{ even}}} \Det(R_i) \right) \circ
\phi_T = \phi_{T'} \circ \left( \bigotimes_{i \;\mbox{\scriptsize{
        odd}}} \Det(R_i) \right).
\end{equation}
This immediately follows from the fact that if $S$ satisfies
(\ref{inve}) for $T$ then the sequence $S'$ defined by $S_i' :=
R_{i-1} S_i R_i^{-1}$ satisfies (\ref{inve}) for $T'$.

The exact sequence $A$
\[
0 \rightarrow A_1 \longrightarrow A_1 \oplus A_2 \longrightarrow A_2
\oplus A_3 \longrightarrow \dots \longrightarrow A_{n-1} \oplus A_n
\longrightarrow A_n \rightarrow 0,
\]
consisting of inclusions and projections produces the isomorphism
$\phi_A$ given by
\[
\phi_A ( \alpha_1 \otimes (\alpha_2 \wedge \alpha_3) \otimes \dots ) =
(\alpha_1 \wedge \alpha_2) \otimes (\alpha_3 \wedge \alpha_4) \otimes
\dots,
\]
where $\alpha_i$ generates $\Det(A_i)$. 
In the case of a general exact sequence $T$, if $\Det(X_i)$ is
generated by ${T_{i-1}}_* \alpha_{i-1} \wedge \alpha_i$, the isomorphism 
$\phi_T$ is given by
\begin{equation}
\label{expli}
\phi_T ( \alpha_1 \otimes ({T_2}_* \alpha_2 \wedge \alpha_3) \otimes
\dots) = ({T_1}_* \alpha_1 \wedge \alpha_2) \otimes 
 ({T_3}_* \alpha_3 \wedge \alpha_4) \otimes \dots \; .
\end{equation}

For every subset $J\subset \{1,\dots, n\}$, the isomorphism $\phi_T$ induces an isomorphism\footnote{The convention here is that a tensor product over an empty set of indices produces the field of scalars, $\R$ or $\C$.}
\[
\phi_T^J : \bigotimes_{j\in J} \Det(X_j)^{*(j-1)} \cong \bigotimes_{j\notin J} \Det(X_j)^{*(j)},
\]
where $V^{*(j)}$ denotes $V$ when $j$ is even, $V^*$ when $j$ is odd. The isomorphism $\phi_T^J$ is defined by tensorizing $\phi_T$ by the identity on $\Det (X_j)^*$, for every even $j\in J$ and every odd $j\notin J$, and identifying each tensor product $V\otimes V^*$ with the field of scalars - $\R$ or $\C$ - by the duality pairing. The original isomorphism $\Phi_T$ corresponds to choosing $J$ to be the subset of $\{1,\dots,n\}$ consisting of odd numbers. 

Naturality now means that if we are given the isomorphism (\ref{isoexa}) between the exact sequences $T$ and $T'$, there holds
\begin{equation}
\label{natural2}
\left( \bigotimes_{j\notin J} \Det(R_j^{(-1)^j})^{*(j)} \right) \circ
\phi_T^J = \phi_{T'}^J \circ \left( \bigotimes_{j\in J} \Det(R_j^{(-1)^{j-1}})^{*(j-1)} \right).
\end{equation}
Since the naturality property only involves pairs of exact sequences where the corresponding vector spaces have identical dimension, we still achieve naturality if we multiply each $\Phi_T^J$ by a non-zero number which depends only on the dimensions of $X_1,\dots,X_n$. In particular, we are interested in changing $\phi_T^J$ just by a factor $\pm 1$. We can summarize the above discussion into the following:

\begin{prop}
\label{natup}
For any choice of the subset $J\subset \{1,\dots,n\}$ and of
the function $\sigma:\N^n \rightarrow \{-1,1\}$, the transformation $T \mapsto \phi_T^{J,\sigma}$ which associates the isomorphism
\[
\phi_T^{J,\sigma} = \sigma(\dim X_1,\dots,\dim X_n) \, \phi_T^J : 
\bigotimes_{j\in J} \Det(X_j)^{*(j-1)} \cong \bigotimes_{j\notin J} \Det(X_j)^{*(j)},
\]
to the exact sequence $T$ of finite dimensional vector spaces
\[
0 \stackrel{T_0}{\longrightarrow} X_1 \stackrel{T_1}{\longrightarrow}
X_2 \stackrel{T_2}{\longrightarrow} \dots
\stackrel{T_{n-1}}{\longrightarrow} X_n
\stackrel{T_n}{\longrightarrow} 0,
\]
is natural with respect to isomorphisms of exact sequences.
\end{prop}

We conclude this section by considering the particular case of 
a sequence with four terms, that is the exact sequence $T$ 
\[
0 \rightarrow \ker T \stackrel{i}{\longrightarrow} 
X \stackrel{T}{\longrightarrow} Y
\stackrel{\pi}{\longrightarrow} \coker T \longrightarrow 0.
\]
induced by some homomorphism $T:X \rightarrow Y$ between finite
dimensional vector spaces. We shall associate to such an exact sequence 
the isomorphism
\begin{equation}
\label{bismut}
\psi_T : \Det(\ker T) \otimes \Det(\coker T)^* \cong \Det(X) \otimes
\Det(Y)^*, \quad \alpha \otimes (\pi_* \gamma)^* \mapsto (i_*(\alpha) \wedge
\beta) \otimes (\gamma \wedge T_* \beta)^*,  
\end{equation}
where $\alpha$ generates $\Det (\ker T)$, $i_*(\alpha) \wedge \beta$
generates $\Det(X)$, $\gamma \wedge T_* \beta$ generates $\Det(Y)$, and the
superscript $*$ indicates the dual generator. Comparing this definition with formula (\ref{expli}), we notice an exchange in the exterior product between the generators $\gamma$ and $T_* \beta$ 
$\psi_T$, which produces a sign change equal to
\[
(-1)^{(\dim \ran T)(\dim \coker T)}.
\]
Therefore, $\psi_T$ coincides with the isomorphism $\phi_T^{J,\sigma}$, where
\[
J=\{1,4\}, \quad \sigma(d_1,d_2,d_3,d_4) = (-1)^{d_4 (d_3 - d_4)}.
\]

Let $X=X_0 \oplus Z$, $Y=Y_0\oplus Z$ be finite dimensional vector
spaces, and let $T\in \mathrm{L}(X,Y)$ be such that $TX_0 \subset
Y_0$ and $QT|_Z$ is an automorphism of $Z$, where $Q:Y_0 \oplus Z
\rightarrow Z$ denotes the projection.  
Denote by $T_0\in \mathrm{L}(X_0,Y_0)$ the restriction of $T$. The fact that
$QT|_Z$ is an isomorphism implies that there are natural identifications
\begin{equation}
\label{idf}
\ker T_0 \cong \ker T, \quad \coker T_0 \cong \coker T.
\end{equation}
Let 
\[
\psi_T : \Det(\ker T) \otimes \Det(\coker T)^* \cong \Det(X) \otimes
\Det(Y)^*
\]
and
\[
\psi_{T_0} : \Det(\ker T_0) \otimes \Det(\coker T_0)^*
\cong \Det(X_0) \otimes \Det(Y_0)^*
\]
be the isomorphism defined above. By the identifications (\ref{idf})
we can consider the composition
\[
\psi_T \circ \psi_{T_0}^{-1} : \Det(X_0) \otimes \Det(Y_0)^* \cong
\Det(X) \otimes \Det(Y)^*.
\]
The next lemma justifies the choice of the sign in the definition of $\psi_T$:

\begin{lem}
\label{restri}
Let $x_0$, $y_0$, $z$ be generators of $\Det(X_0)$, $\Det(Y_0)$,
$\Det(Z)$, respectively. Then
\[
\psi_T \circ \psi_{T_0}^{-1} (x_0 \otimes y_0^*) = 
\det (QT|_Z) \, (x_0 \wedge z) \otimes (y_0\wedge z)^*.
\]
In particular, the function mapping a homomorphism $T\in \mathrm{L}(X,Y)$ with the above properties into $\psi_T \circ \psi_{T_0}^{-1}$ depends analytically on $T$.
\end{lem}

\begin{proof}
We can write $x_0=\alpha \wedge x_1$, where $\alpha$ is a generator of
$\Det(\ker T_0)$, and $y_0=y_1 \wedge {T_0}_* x_1$, where $y_1$ is a
generator of $\coker T_0\cong \coker T$, identified with a subspace of 
$Y_0$. Then
\begin{eqnarray*}
\psi_{T_0} (\alpha \otimes y_1^*) = (\alpha \wedge x_1) \otimes
(y_1 \wedge {T_0}_* x_1)^* = x_0 \otimes y_0^*, \\
\psi _T(\alpha \otimes y_1^*) = (\alpha \wedge x_1 \wedge z) \otimes 
(y_1 \wedge T_*(x_1 \wedge z))^* = (x_0 \wedge z) \otimes 
(y_1 \wedge T_*(x_1 \wedge z))^*.
\end{eqnarray*}
Notice that
\begin{eqnarray*}
y_1 \wedge T_*(x_1 \wedge z) = y_1 \wedge {T_0}_* x_1 \wedge T_* z =
y_0 \wedge QT_* z + y_0 \wedge
(I-Q)T_* z = \det (QT|_Z) y_0 \wedge z.
\end{eqnarray*}
It follows that
\[
\psi_T(\alpha \otimes y_1^*) = \det
  (QT|_Z) (x_0 \wedge z) \otimes (y_0 \wedge z)^*,
\]
so
\[
\psi_T \circ \psi_{T_0}^{-1} (x_0 \otimes y_0^*) = 
\det (QT|_Z) \, (x_0 \wedge z) \otimes (y_0\wedge z)^*.
\]
as claimed.
\end{proof} \qed

Another useful property of the isomorphism $\psi_T$ defined in
(\ref{bismut}) is that
\begin{equation}
\label{adj}
\psi_{T^*}^{-1} = \psi_T^*,
\end{equation}
whereas the corresponding formula relating a general isomorphism $\phi_T^{J,\sigma}$ to $\phi_{T^*}^{J,\sigma}$ is more complicated, and may involve a change of sign.

\section{The determinant bundle over the space of Fredholm pairs}
\label{detbun}

The determinant bundle over $\mathrm{Fp}(H)$ is an analytic line bundle
with base space $\mathrm{Fp}(H)$, and fiber at $(V,W)$
\[
\Det (V,W) := \Det (V\cap W) \otimes \Det\left( \frac{H}{V+W}
\right)^*.
\]
It is denoted by
\[
p:\Det (\mathrm{Fp}(H)) \rightarrow \mathrm{Fp}(H), \quad \Det(\mathrm{Fp}(H)) :=
\bigsqcup_{(V,W) \in \mathrm{Fp}(H)} \Det (V,W).
\]
The construction of an analytic structure on the total space
$\Det(\mathrm{Fp}(H))$, making $p$ an analytic bundle map, is not immediate
because the subspaces $V\cap W$ and $V+W$ do no depend even
continuously on $(V,W)\in \mathrm{Fp}(H)$. We will describe such an analytic
structure first on the counter image by $p$ of those components of
$\mathrm{Fp}(H)$ where one of the two spaces is finite dimensional, and then on
the other components.

\paragraph{Analytic structure on $\mathbf{p^{-1}(\mathrm{Fp}(H)\setminus \mathrm{Fp}^*(H))}$.}
If $(V,W)\in \mathrm{Gr}_{n,\infty}(H)\times \mathrm{Gr}_{\infty,m}(H)$, $n,m\in \N$, 
the exact sequence of inclusions and quotient projections
\[
0 \rightarrow V\cap W \rightarrow V \rightarrow \frac{H}{W}
\rightarrow \frac{H}{V+W} \rightarrow 0
\]
yields to the isomorphism
\begin{equation}
\label{duno}
\psi: \Det(V,W) \cong \Det (V) \otimes \Det \left( \frac{H}{W} \right)^*,
\end{equation}
defined in (\ref{bismut}). The lines $\Det(V) \otimes \Det(H/W)^*$ are 
clearly the fibers of an analytic bundle over $\mathrm{Gr}_{n,\infty}(H)\times 
\mathrm{Gr}_{\infty,m}(H)$, and we can endow $p^{-1}(\mathrm{Gr}_{n,\infty}(H)\times 
\mathrm{Gr}_{\infty,m}(H))$ with the analytic structure which makes
(\ref{duno}) an isomorphism of analytic bundles. The same construction 
works for $p^{-1}(\mathrm{Gr}_{\infty,n}(H)\times \mathrm{Gr}_{m,\infty}(H))$. 

\paragraph{Analytic structure on $\mathbf{p^{-1}(\mathrm{Fp}^*(H))}$.}
We proceed to define a structure of analytic bundle on the restriction
of $p$ to $p^{-1}(\mathrm{Fp}^*(H))$. Consider the covering of $\mathrm{Fp}^*(H)$
consisting of the open sets
\[
\mathcal{U}_X := \set{(V,W)\in \mathrm{Fp}^*(H)}{X+V+W=H, \; X\cap V=(0)},
\]
for $X$ varying among all finite dimensional linear subspaces of
$H$. If $X$ is finite dimensional, the map $V\mapsto X+V$ is analytic 
on the space of $V\in \mathrm{Gr}(H)$ such that $X\cap V=(0)$: indeed in this
case $X+V$ is closed, so 
$X+V = (X^{\perp} \cap V^{\perp} )^{\perp}$ and $X^{\perp} + V^{\perp}
= (X\cap V)^{\perp} = H$, hence the analyticity of the
sum follows from the analyticity of the transverse intersection. As a
consequence,  the map $\mathcal{U}_X \rightarrow \mathrm{Gr}(H)$, $(V,W)\mapsto
(X+V)\cap W$, is analytic. Therefore, we can consider the analytic line bundle
$p_X : \mathcal{E}_X \rightarrow \mathcal{U}_X$, whose fiber at
$(V,W)$ is the line
\[
p_X^{-1} (\{(V,W)\}) := \Det ((X+V)\cap W) \otimes \Det (X)^*.
\]
Let $(V,W)\in \mathcal{U}_X$. Since $X\cap V=(0)$, we can define a
linear map $(X+V)\cap W \rightarrow X$ to be the following
composition of inclusions, quotient projections, and their inverses:
\[
(X+V)\cap W \rightarrow X+V \rightarrow
\frac{X+V}{V} \stackrel{\cong}{\rightarrow} X.
\]
The kernel of this linear map is $V\cap W$, while its range is
$X\cap (V+W)$, so using also the identity $X+V+W=H$, we obtain the
exact sequence
\[
0 \rightarrow V\cap W \rightarrow (X+V)\cap W
\longrightarrow X \rightarrow
\frac{H}{V+W} \rightarrow 0.
\]
Therefore formula (\ref{bismut}) produces the isomorphism
\[
\psi_{(V,W)}^X: \Det(V,W) \rightarrow \Det((X+V)\cap W) \otimes
\Det(X)^*,
\]
and we can endow $p^{-1}(\mathcal{U}_X)$ with the analytic structure
which makes the map $p^{-1}(\mathcal{U}_X) \rightarrow \mathcal{E}_X$,
$\xi \mapsto \psi_{(V,W)}^X \xi$ for $p(\xi)=(V,W)$, an isomorphism
of analytic bundles. 

In order to have an analytic structure on the whole $p^{-1}(\mathrm{Fp}^*(H))$,
there remains to show that if $X_1,X_2$ are finite dimensional linear
subspaces of $H$, the transition map
\[
\psi_{X_1}^{X_2} : p_{X_1}^{-1} (\mathcal{U}_{X_1} \cap
\mathcal{U}_{X_2}) \rightarrow p_{X_2}^{-1} (\mathcal{U}_{X_1} \cap
\mathcal{U}_{X_2}), \quad \xi \mapsto \psi_{(V,W)}^{X_2}
(\psi_{(V,W)}^{X_1})^{-1} \xi, \quad \mbox{for } p_{X_1}(\xi) =
(V,W),
\]
is analytic. When $X_1\subset X_2$ or $X_2\subset X_1$, this fact
follows immediately from Lemma \ref{restri}. In the general case, let
$(V,W)\in \mathcal{U}_{X_1} \cap \mathcal{U}_{X_2}$. Then there exists
a finite dimensional linear subspace $X\subset H$ such that
\begin{equation}
\label{ddue}
X+V+W=H, \quad (X_1+X)\cap V = (0), \quad (X_2+X) \cap V=(0).
\end{equation}
Indeed, setting $k=\codim (V+W)$ and recalling that $V$ has infinite
codimension, we have that the set of $X\in \mathrm{Gr}_{k,\infty}(H)$
satisfying each of the identities (\ref{ddue}) is open and
dense. Therefore,
\[
(V,W) \in \mathcal{U}_{X_1} \cap \mathcal{U}_{X_1+X} \cap
\mathcal{U}_{X} \cap \mathcal{U}_{X_2+X} \cap \mathcal{U}_{X_2},
\]
and on the inverse image by $p$ of 
such a neighborhood $\psi_{X_1}^{X_2}$ is the composition of
analytic maps, by the case seen above. Hence $\psi_{X_1}^{X_2}$ is
analytic on $p^{-1}(\mathcal{U}_{X_1} \cap \mathcal{U}_{X_2})$.

\paragraph{Transposition.} 
Note that, although the fiber $\Det (V,W)$ depends symmetrically on
$(V,W)$, in our construction of the analytic structure on
$p^{-1}(\mathrm{Fp}^*(H))$ such a symmetry is lost. However, if we exchange the
role of $V$ and $W$ in the definition above, we obtain the same
analytic structure. In other words, the map $\Det(\mathrm{Fp}(H))\rightarrow
\Det(\mathrm{Fp}(H))$ which lifts the transposition involution
\[
\tau: \mathrm{Fp}(H) \rightarrow \mathrm{Fp}(H), \quad
(V,W)\mapsto (W,V),
\]
and is the identity on the fibers, is an analytic bundle isomorphism. This fact is obvious for the components where one of the spaces is finite dimensional, so we just have to check it locally at some $(V_0,W_0)\in \mathrm{Fp}^*(H)$. We can choose a finite dimensional linear subspace $X\subset H$ such that $X\cap (V_0+W_0)=(0)$ and $X+V_0+W_0=H$. Then both
$(V_0,W_0)$ and $(W_0,V_0)$ belong to the $\tau$-invariant open set
\[
\mathcal{U}_X \cap \tau(\mathcal{U}_X) =
\set{(V,W)\in \mathrm{Fp}^*(H)}{X\cap V = X \cap W =(0), \; X+V+W=H},
\]
and we must check that the isomorphism
\[
\psi_{(W,V)}^X \circ \bigl(\psi_{(V,W)}^X\bigr)^{-1} : \Det ((X+V)\cap W) \otimes \Det(X)^* \rightarrow \Det ((X+W)\cap V) \otimes \Det(X)^*
\]
depends analytically on $(V,W)$ in such a set. If $(V,W)$ belong to this set, we can find an isomorphism
\[
T_{(V,W)} : (X+W)\cap V \stackrel{\cong}{\longrightarrow} (X+V)\cap W,
\]
that makes the following diagram
\begin{equation*}
\xymatrix{
0 \ar[r] & V\cap W \ar[r] \ar[d]_{\mathrm{id}} & (X+V)\cap W \ar[r] \ar[d]_{T_{(V,W)}}& X \ar[r]  \ar[d]^{\mathrm{id}} & \frac{H}{V+W} \ar[r] \ar[d]^{\mathrm{id}} & 0 \\ 
0  \ar[r] & W\cap V \ar[r]  & (X+W)\cap V \ar[r]  & X \ar[r]  & \frac{H}{W+V} \ar[r]  & 0 }
\end{equation*}
commutative. Such an isomorphism is uniquely determined by the choice of two linear complements of $V\cap W$ in $(X+V)\cap W$ and $(X+W)\cap V$, and if we choose these to be the orthogonal complements, $T_{(V,W)}$ depends analytically on $(V,W)$.  By the naturality property stated in Proposition \ref{natup}, we have
\[
\psi_{(W,V)}^X \circ \bigl(\psi_{(V,W)}^X\bigr)^{-1} = \Det(T_{(V,W)}) \otimes \mathrm{id},
\]
which depends analytically on $(V,W)\in \mathcal{U}_X \cap \tau(\mathcal{U}_X)$, proving the claim.
 
This fact shows that it is possible to define a determinant bundle also on the space of unordered Fredholm pairs.

\paragraph{Nontriviality.} 
The determinant bundle $p:\Det(\mathrm{Fp}(H))
\rightarrow \mathrm{Fp}(H)$ is nontrivial on every connected component of
$\mathrm{Fp}(H)$, with the exception of the trivial components
$\mathrm{Gr}_{0,\infty}(H) \times \mathrm{Gr}_{\infty,0}(H)$ and $\mathrm{Gr}_{\infty,0}(H)
\times \mathrm{Gr}_{0,\infty}(H)$. We shall check this fact for the connected
components 
$\mathrm{Fp}_k^*(H)$, $k\in \Z$, the case of the other components being
analogous. First notice that if $H_0$ is 2-dimensional, then the
determinant bundle on $\mathrm{Gr}_{1,1}(H_0)$ is nontrivial, being isomorphic
to the tautological line bundle on the projective line. Fix a
splitting $H=H_0 \oplus H_1$, $\dim H_0 = 2$, a Fredholm pair
$(V,W)\in \mathrm{Fp}^*_{k-1}(H_1)$ and a generator $\xi$ of $\Det(V \cap
W) \otimes \Det(H_1/(V+W))^*$, so that for every $L\in \mathrm{Gr}_{1,1}(H_0)$,
\[
(L\oplus V) \cap (H_0 \oplus W) = L \oplus (V\cap W), \quad
\frac{H}{(L\oplus V) + (H_0 \oplus W)} = \frac{H_1}{V+W}.
\]
Then the embedding
\[
f:\mathrm{Gr}_{1,1}(H_0) \rightarrow \mathrm{Fp}_k(H), \quad L \mapsto (L \oplus V, H_0
\oplus W) 
\]
can be lifted to an analytic bundle isomorphism
\[
F: \Det(\mathrm{Gr}_{1,1}(H_0)) \rightarrow p^{-1}(f(\mathrm{Gr}_{1,1}(H_0))) \subset
\Det(\mathrm{Fp}_k^*(H)), \quad\quad \eta \mapsto \eta \otimes \xi,
\]
and the determinant bundle is nontrivial on $\mathrm{Fp}^*_k(H)$. 

\paragraph{The determinant bundle over $\mathbf{Gr_c(V,H)}$.} 
Let $V\in Gr_{\infty,\infty}(H)$. Then the Grassmannian $\mathrm{Gr}_c(V,H)$ of compact perturbations of $V$ has a natural inclusion into $\mathrm{Fp}(H)$ given by
\begin{equation}
\label{incl}
\mathrm{Gr}_c(V,H) \hookrightarrow \mathrm{Fp}(H), \quad W \mapsto (W,V^{\perp}).
\end{equation}
The above map allows us to consider the pull-back of the determinant
bundle $\Det(\mathrm{Fp}(H))\rightarrow \mathrm{Fp}(H)$ on $\mathrm{Gr_c}(V,H)$, obtaining an analytic  line
bundle $\Det(\mathrm{Gr_c}(V,H))\rightarrow \mathrm{Gr_c}(V,H)$ with fiber at $W\in \mathrm{Gr_c}(V,H)$
\[
\Det(W\cap V^{\perp}) \otimes \Det\left(
\frac{H}{W + V^{\perp}}\right)^* \cong \Det
\bigl( (W\cap V^{\perp}) \oplus (W^{\perp} \cap V) \bigr).
\]
The argument for the nontriviality used above shows
that the above line bundle is nontrivial on every connected component
$\mathrm{Gr}_{\mathrm{c},n}(V,H)$, for $n\in \Z$. 

\section{The determinant bundle over the space of Fredholm operators.} 

Let $\mathrm{Fr}(H_1,H_2)$ be the open subset of $\mathrm{L}(H_1,H_2)$
consisting of Fredholm operators. The inclusion 
\[
i: \mathrm{Fr}(H_1,H_2) \rightarrow \mathrm{Fp}(H_1 \times H_2), \quad T \mapsto (\graf
T, H_1 \times (0)),
\]
allows to define the determinant bundle over $\mathrm{Fr}(H_1,H_2)$ as the
pull-back of the determinant bundle over $\mathrm{Fp}(H_1\times H_2)$ by the
map $i$. We obtain a line bundle
\[
q:\Det(\mathrm{Fr}(H_1,H_2)) \rightarrow \mathrm{Fr}(H_1,H_2),
\]
whose fiber at $T$ is
\[
\Det(T) := \Det(\ker T) \otimes \Det(\coker T)^*.
\]
This is the determinant bundle over the space of Fredholm operators
defined by Quillen in \cite{qui85} (see also \cite{bf86}). 
A direct construction of the
analytic bundle structure on $\Det(\mathrm{Fr}(H_1,H_2))$ goes as follows. Let
$X$ be a finite dimensional linear subspace of $H_2$, and consider the
open set
\begin{equation}
\label{fredcov}
\mathcal{U}_X(H_1,H_2) := \set{T\in \mathrm{Fr}(H_1,H_2)}{T \mbox{ is transverse to } X}.
\end{equation}
If $T\in \mathcal{U}_X(H_1,H_2)$, the exact sequence 
\[
0 \rightarrow \ker T \longrightarrow T^{-1} X
\stackrel{T}{\longrightarrow} X \longrightarrow \coker T \rightarrow 0
\]
induces an isomorphism
\begin{equation}
\label{atlfred}
\psi_T : \Det(T) \rightarrow \Det(T^{-1} X) \otimes \Det(X)^*
\end{equation}
by formula (\ref{bismut}). Since the lines $\Det(T^{-1} X) \otimes \Det(X)$
are the fibers of an analytic line bundle over
$\mathcal{U}_X(H_1,H_2)$, we obtain an analytic line bundle structure
for $\Det(\mathrm{Fr}(H_1,H_2))$ over $\mathcal{U}_X(H_1,H_2)$. Lemma
\ref{restri} can again be used to show that the transition maps are
analytic. 
 
An alternative construction for the analytic structure of $q$ is the following.
If $T$ is in $\mathrm{Fr}(H_1,H_2)$, the operator $T^*T\in \mathrm{L}(H_1)$ is
  self-adjoint, positive, and Fredholm. In particular, if $\epsilon>0$
  is small enough, the spectrum of $T^*T$ has finite multiplicity in
  $[0,\epsilon]$, meaning that $\sigma(T^*T) \cap [0,\epsilon]$ is a
  finite set of eigenvalues with finite multiplicity. Denote by
  $V_{\epsilon}(T^*T)$ the corresponding finite dimensional
  eigenspace. Therefore the sets
\[
\mathcal{V}_{\epsilon}(H_1,H_2) := \set{T\in \mathrm{Fr}(H_1,H_2)}{\epsilon\notin
  \sigma(T^*T) \mbox{ and } \sigma(T^*T)\cap [0,\epsilon] \mbox{
  has finite multiplicity}}, \quad \epsilon\geq 0,
\]
constitute an open covering of $\mathrm{Fr}(H_1,H_2)$. Clearly, the operators 
$T^*T|_{(\ker T)^{\perp}} : (\ker T)^{\perp} \rightarrow (\ker
T)^{\perp} = \ran T^*$, and $TT^*|_{(\ker T^*)^{\perp}}  : (\ker
T^*)^{\perp} \rightarrow (\ker T^*)^{\perp} = \ran T$, are conjugated
by $T|_{(\ker T)^{\perp}} :  (\ker T)^{\perp} \rightarrow (\ker
T^*)^{\perp}$. Since $\ker T^*T=\ker T$ and $\ker TT^*=\ker T^*$, it follows
that any $\lambda>0$ belongs to $\sigma(T^*T)$ if and only if it
belongs to $\sigma(TT^*)$, and that it is
an eigenvalue of $T^*T$ corresponding to the
eigenvector $x$ if and only if it is an eigenvalue of $TT^*$
corresponding to the eigenvector $Tx$, with the same
multiplicity. Hence, if $T\in \mathcal{V}_{\epsilon}(H_1,H_2)$, 
then $T^*\in \mathcal{V}_{\epsilon}(H_2,H_1)$, and the sequence
\begin{equation}
\label{exseq}
0 \rightarrow \ker T \rightarrow V_{\epsilon}(T^*T)
\stackrel{T}{\longrightarrow} V_{\epsilon}(TT^*) \rightarrow \ker T^*
\rightarrow 0
\end{equation}
is exact. From the isomorphism
\[
\psi_T^{\epsilon}: \Det (T) \cong \Det (V_{\epsilon}(T^*T))^* \otimes
\Det(V_{\epsilon}(TT^*))
\]
defined by (\ref{bismut})
and from the analyticity of the maps $T\mapsto V_{\epsilon}(T^*T)$ and
$T \mapsto V_{\epsilon}(TT^*)$ on $\mathcal{V}_{\epsilon}(H_1,H_2)$, we obtain
the analytic structure of $q$ over $\mathcal{V}_{\epsilon}(H_1,H_2)$. The
analyticity of the transition maps is an immediate consequence of Lemma
\ref{restri}, because if $0\leq \epsilon<\epsilon^{\prime}$ and $T\in
\mathcal{V}_{\epsilon}(H_1,H_2) \cap \mathcal{V}_{\epsilon^{\prime}} 
(H_1,H_2)$, there holds $V_{\epsilon}(T^*T) \subset 
V_{\epsilon^{\prime}}(T^*T)$, and
$V_{\epsilon}(TT^*) \subset V_{\epsilon^{\prime}}(TT^*)$.  

We just mention that an analogous spectral approach provides an
alternative construction of the determinant bundle on $\mathrm{Fp}(H)$.

\paragraph{Left and right action of $\mathbf{\mathrm{GL}(H)}$.} The group $\mathrm{GL}(H_2)$ acts on
$\mathrm{Fr}(H_1,H_2)$ by left multiplication $(G,T)\mapsto GT$, while the
group $\mathrm{GL}(H_1)$  acts on $\mathrm{Fr}(H_1,H_2)$ by right multiplication
$(G,T) \mapsto TG$. Since
\[
\ker GT = \ker T, \quad \coker GT = \tilde{G} \, \coker T,
\]
where $\tilde{G}:H_2/\ran T \rightarrow H_2/G\, \ran T$ is induced by $G$,
the left action lifts to a bundle action on $\Det(\mathrm{Fr}(H_1,H_2))$ defined fiberwise by
\[
\Det (T) \rightarrow \Det(GT) , \quad \xi \otimes \eta^* \mapsto
\xi \otimes \bigl(\Det(\tilde{G}^{-1})^* \eta^*\bigr).
\]
This action is analytic because of the naturality property of Proposition \ref{natup}. Similarly, since
\[
\ker TG = G^{-1} \ker T, \quad \coker TG = \coker T,
\]
the right action lifts to an analytic bundle action on
$\Det(\mathrm{Fr}(H_1,H_2))$ defined fiberwise by
\[
\Det (T) \rightarrow \Det(TG) , \quad \xi \otimes \eta^* \mapsto
\bigl(\Det(G^{-1})\xi \bigr) \otimes \eta^*.
\]

\paragraph{Adjoint.}  
The adjoint map
\[
\mathrm{Fr}(H_1,H_2) \rightarrow \mathrm{Fr}(H_2,H_1), \quad T \mapsto T^*,
\]
has analytic lift to the determinant bundles. Indeed, by using the
identifications
\[
\ker T^* = (\ran T)^{\perp} \cong (\coker T)^*, \quad
\coker T^* = \frac{H_1^*}{(\ker T)^{\perp}} \cong (\ker T)^*,
\]
the lift is defined fiberwise by 
\[
\Det(T) \rightarrow \Det(T^*) , \quad \xi \otimes \eta^* \mapsto \eta^*
\otimes \xi.
\]
Let us check that this lift is analyitic.
If $T$ belongs to the open set $\mathcal{V}_{\epsilon}(H_1,H_2)$, then
$T^*$ belongs to the open set $\mathcal{V}_{\epsilon}(H_2,H_1)$.
The adjoint of the exact sequence (\ref{exseq}) is 
\[
0 \rightarrow \ker T^* \rightarrow V_{\epsilon}(TT^*)
\stackrel{T}{\longrightarrow} V_{\epsilon}(T^*T) \rightarrow \ker T
\rightarrow 0,
\]
which is precisely the sequence producing the isomorphism
\[
\psi_{T^*}^{\epsilon}: \Det(T^*) \cong \Det (V_{\epsilon}(TT^*)) \otimes
\Det(V_{\epsilon}(T^*T))^*.
\]
Then property (\ref{adj}) allows to conclude.

\paragraph{Two analytic sections.} The determinant bundle
$\Det(\mathrm{Fr}(H_1,H_2))$ has the global analytic section
\[
s(T) = \left\{ \begin{array}{ll} 0 & \mbox{if $T$ is not invertible} \\ 1 \in
  \Det(T) = \C \otimes \C^*  & \mbox{if $T$ is
  invertible}. \end{array} \right.
\]
This section vanishes precisely at the non-invertible elements of
$\mathrm{Fr}(H_1,H_2)$. 

The restriction of the determinant bundle $\Det(\mathrm{Fr}(H,H))$ to the space
of self-adjoint Fredholm operators on $H$ has a nowhere vanishing analytic
section $\sigma$ which associates to each self-adjoint Fredholm
operator $T$ the element 
\[
1 \in \Det(\ker T) \otimes \Det (\ker T)^* = \Det(\ker T) \otimes
\Det (\coker T)^* = \Det (T).
\]
In particular, the restriction of the determinant bundle over the
space of self-adjoint Fredholm operators on $H$ is trivial (as observed
by Furutani in \cite{fur04}, Theorem 4.1). 

\paragraph{Composition.} Let us show how the composition of Fredholm
operators
\[
\mathrm{Fr}(H_1,H_2) \times \mathrm{Fr}(H_2,H_3) \rightarrow \mathrm{Fr}(H_1,H_3), \quad 
(S,T) \mapsto TS,
\]
lifts to the determinant bundles, producing the bundle morphism
\[
\Det(\mathrm{Fr}(H_1,H_2)) \otimes \Det(\mathrm{Fr}(H_2,H_3)) \rightarrow
\Det(\mathrm{Fr}(H_1,H_3)),
\]
where the domain is seen as a line bundle over the product
$\mathrm{Fr}(H_1,H_2) \times \mathrm{Fr}(H_2,H_3)$. 

If $S\in \mathrm{Fr}(H_1,H_2)$ and $T\in \mathrm{Fr}(H_2,H_3)$, we have the exact
sequence
\begin{equation}
\label{exa}
0 \rightarrow \ker S \longrightarrow \ker TS
\stackrel{S}{\longrightarrow} \ker T \stackrel{\pi}{\longrightarrow}
\coker S \stackrel{T}{\longrightarrow} \coker TS \longrightarrow
\coker T \rightarrow 0,
\end{equation}
where $\pi$ denotes the restriction of the quotient projection. This
exact sequence induces the isomorphism
\begin{equation}
\label{compo}
\begin{aligned}
\phi_{S,T} : \Det(S) \otimes \Det(T)  & \rightarrow  \Det(TS), \\ 
\alpha_1 \otimes (\pi_* \alpha_3 \wedge \alpha_4)^* \otimes (\alpha_3
\wedge S_* \alpha_2) \otimes \alpha_5^*  & \mapsto  (-1)^{\rho(S,T)} 
(\alpha_1 \wedge
\alpha_2) \otimes (\alpha_5 \otimes T_* \alpha_4)^*, \end{aligned} 
\end{equation}
with
\[
\rho(S,T)   := (\dim \ker S +
  \dim \ker TS)(\dim \coker T + \dim \coker TS), 
\]
where $\alpha_1$ generates $\Det(\ker S)$, $\alpha_1\wedge \alpha_2$
generates $\Det (\ker TS)$, $\alpha_3 \wedge S_* \alpha_2$ generates
$\Det(\ker T)$, $\pi_* \alpha_3 \wedge \alpha_4$ generates $\Det
(\coker S)$, $\alpha_5 \wedge T_* \alpha_4$ generates
$\Det(\coker TS)$, and $\alpha_5$ generates $\Det (\coker T)$. If we denote by $(S,T)$ the exact sequence (\ref{exa}), the isomorphism $\phi_{S,T}$ coincides with the natural isomorphism $\phi_{(S,T)}^{J,\sigma}$ of section \ref{Dets}, with
\begin{equation}
\label{sgnch1}
J = \{1,3,4,6\}, \quad \sigma(d_1,d_2,d_3,d_4,d_5,d_6) = (-1)^{d_4(d_2-d_1) + d_6 (d_5 - d_6)}.
\end{equation}

Let us check that the bundle map defined by the isomorphisms
$\phi_{S,T}$ is analytic. It is convenient to see the analytic structure on the determinant bundle of $\mathrm{Fr}(H_1,H_2)$ in terms of the open sets $\mathcal{U}_X(H_1,H_2)$ and of the isomorphisms $\psi_T$ introduced in (\ref{fredcov}) and (\ref{atlfred}).

Let $Y$ be a finite dimensional linear
subspace of $H_3$. Notice that the composition $TS$ is transverse to
$Y$ if and only if $T$ is transverse to $Y$ and $S$ is transverse to
$X:= T^{-1} Y$. In particular, the composition maps
$\mathcal{U}_X(H_1,H_2) \times \mathcal{U}_Y(H_2,H_3)$ into
$\mathcal{U}_Y(H_1,H_3)$. Let 
\begin{eqnarray*}
\psi_S : \Det(S) & \rightarrow & \Det(S^{-1} X) \otimes \Det(X)^*, \\
\psi_T : \Det(T) & \rightarrow & \Det(T^{-1} Y) \otimes \Det(Y)^* = \Det(X) \otimes \Det(Y)^*, \\
\psi_{TS} : \Det(TS) & \rightarrow & \Det((TS)^{-1} Y) \otimes \Det(Y)^* = \Det(S^{-1} X) \otimes \Det(Y)^*, 
\end{eqnarray*}
be the isomorphisms defined in (\ref{atlfred}).
We must show that the upper horizontal
isomorphism which makes the following diagram commutative,
\begin{equation*}
\begin{CD}
\Det(S^{-1}X) \otimes \Det(X)^* \otimes \Det(X) \otimes \Det(Y)^* @>>>
\Det(S^{-1}X) \otimes \Det(Y)^* \\ @A{\psi_S \otimes \psi_T}AA
@AA{\psi_{TS}}A \\ \Det(\ker S) \otimes \Det (\coker S)^* \otimes
\Det(\ker T) \otimes \Det (\coker T)^* @>{\phi_{S,T}}>> \Det(\ker TS)
\otimes \Det (\coker TS)^*
\end{CD}\end{equation*}
depends analytically on $(S,T)\in  \mathcal{U}_X(H_1,H_2) \times
\mathcal{U}_Y(H_2,H_3)$.

Let $\alpha_1,\alpha_2,\alpha_3,\alpha_4,\alpha_5$ be as above. Let
$\zeta$ be a generator of $\Det(Y)$. Since $TS$ is transverse to $Y$,
we can find an element $\alpha_0\in
\Lambda_*((TS)^{-1}Y)=\Lambda_*(S^{-1}X)$ such that $\zeta = \alpha_5
\wedge T_* \alpha_4 \wedge T_*S_* \alpha_0$. Then $\xi :=
\alpha_1 \wedge \alpha_2 \wedge \alpha_0$ generates
$\Det(S^{-1}X)$. Moreover, 
$\eta := \alpha_3
\wedge \alpha_4 \wedge S_* \alpha_2 \wedge S_* \alpha_0$ generates
$\Det(X)$.  By (\ref{bismut}), we have
\begin{eqnarray*}
\psi_S ( \alpha_1 \otimes (\pi_* \alpha_3 \wedge \alpha_4) ^*) & = &
(\alpha_1 \wedge \alpha_2 \wedge \alpha_0) \otimes (\alpha_3 \wedge
\alpha_4 \wedge S_* \alpha_2 \wedge S_* \alpha_0)^* = \xi \otimes
\eta^*, \\
\psi_T ( (\alpha_3 \wedge S_* \alpha_2) \otimes \alpha_5^*) & = &
(\alpha_3 \wedge S_* \alpha_2 \wedge \alpha_4 \wedge S_* \alpha_0)
\otimes (\alpha_5 \wedge T_* \alpha_4 \wedge T_*S_* \alpha_0)^* \\ & = &
(-1)^{|\alpha_2| \, |\alpha_4|} \eta \otimes \zeta^*, \\
\psi_{TS} ( (\alpha_1 \wedge \alpha_2) \otimes (\alpha_5 \wedge T_*
\alpha_4)^*) & = &(\alpha_1 \wedge \alpha_2 \wedge \alpha_0) \otimes
(\alpha_5 \wedge T_* \alpha_4 \wedge T_*S_* \alpha_0)^* = \xi
\otimes \zeta^*.
\end{eqnarray*}
Since the parity of $|\alpha_2| \, |\alpha_4|$ equals the parity of
$\rho(S,T)$, the above formulas imply that
\[
\psi_{TS} \circ \phi_{S,T} \circ (\psi_S \otimes \psi_T)^{-1}
(\xi \otimes \eta^* \otimes \eta \otimes \zeta^*) =
\xi \otimes \zeta^*,
\]
so this isomorphism depends analytically on $(S,T)$, as claimed. 

\paragraph{Associativity.} The commutative diagram
\begin{equation*}
\begin{CD}
\mathrm{Fr}(H_1,H_2) \times \mathrm{Fr}(H_2,H_3) \times \mathrm{Fr}(H_3,H_4) @>>> \mathrm{Fr}(H_1,H_2)
\times \mathrm{Fr}(H_2,H_4)  \\ 
@VVV @VVV  \\
\mathrm{Fr}(H_1,H_3) \times \mathrm{Fr}(H_3,H_4) @>>> \mathrm{Fr}(H_1,H_4) 
\end{CD}
\end{equation*}
lifts to a commutative diagram between the corresponding determinant
bundles. In other words, if $T_i \in \mathrm{Fr}(H_i,H_{i+1})$ then
the diagram 
\begin{equation}
\label{diago}
\begin{CD}
\Det(T_1) \otimes \Det(T_2) \otimes \Det(T_3) @>{\mathrm{id}\otimes
  \phi_{T_2,T_3}}>> \Det(T_1) \otimes \Det(T_3 T_2) \\ @VV{\phi_{T_1,T_2} \otimes
  \mathrm{id}}V @VV{\phi_{T_1,T_3 T_2}}V \\ \Det(T_2 T_1) \otimes \Det(T_3)
  @>{\phi_{T_2 T_1,T_3}}>> \Det(T_3 T_2 T_1)
\end{CD}
\end{equation} 
commutes. Although completely elementary, the proof of the commutativity of the above diagram is quite long. It is contained in  Appendix B.

\section{Operations on the determinant bundle over Fredholm pairs}
\label{odfp}

Besides the transposition involution already considered in section \ref{detbun}, other operations on the space of Fredholm pairs can be lifted to the determinant bundle.

\paragraph{The action of $\mathbf{\mathrm{GL}(H)}$.} 
The analytic action of $\mathrm{GL}(H)$ on
$\mathrm{Fp}(H)$, $(T,(V,W))\mapsto (TV,TW)$, lifts to an analytic bundle
action of $\mathrm{GL}(H)$ on $\Det(\mathrm{Fp}(H))$, defined fiberwise by
\[
\Det(V,W) \rightarrow \Det(TV,TW), \quad \xi\otimes \eta^* \mapsto
\Det(T) \xi \otimes \Det(\tilde{T}^{-1})^* \eta^*,
\]
for $\xi\in \Det(V\cap W)$, $\eta^* \in \Det(H/(V+W))^*$, where $\tilde{T}\in
\mathrm{L}(H/(V+W),H/(TV+TW))$ is induced by the linear operator
$T\in \mathrm{GL}(H)$. The analyticity of such an action follows immediately from
the naturality property stated in Proposition \ref{natup}.

\paragraph{Sum.}
If $X,Y,Z$ are finite dimensional linear spaces, the isomorphism 
\[
S(X,Y):\Det(X) \otimes \Det(Y) \rightarrow \Det(X \oplus Y), \quad \xi \otimes \eta \mapsto \xi \wedge \eta,
\]
is induced by the exact sequence
\[
0 \rightarrow X \rightarrow X \oplus Y \rightarrow Y \rightarrow 0.
\]
It is readily seen that the diagram
\begin{equation}
\label{asc}
\begin{CD}
\Det(X) \otimes \Det(Y) \otimes \Det(Z) @>{S(X,Y)\otimes
  \mathrm{id}_{\Det(Z)}}>> \Det(X \oplus Y) \otimes \Det(Z) \\ 
@V{\mathrm{id}_{\Det(X)}\otimes
  S(Y,Z)}VV @VV{S(X\oplus Y,Z)}V \\ \Det(X) \otimes \Det(Y\oplus Z) @>{S(X,Y\oplus Z)}>>
  \Det(X\oplus Y \oplus Z) 
\end{CD}
\end{equation}
commutes. We would like to extend
this construction to Fredholm pairs. 
 
Let $X\in \mathrm{Gr}(H)$ be finite dimensional, and let $(V,W)\in \mathrm{Fp}(H)$ be such that $X\cap V=(0)$. The isomorphism
\[
S(X,(V,W)) : \Det(X) \otimes \Det(V,W) \rightarrow \Det(X+V,W).
\]
is induced by the exact sequence
\begin{equation}
\label{es1}
0 \rightarrow V \cap W \rightarrow (X+V)\cap W
\longrightarrow \frac{X+V}{V} \cong X \rightarrow
\frac{H}{V+W} \rightarrow \frac{H}{X+V+W} \rightarrow 0.
\end{equation}
More precisely, if $(X,V,W)$ denotes the above exact sequence, $S(X,(V,W))$ is defined to be the isomorphism $\phi_{(X,V,W)}^{J,\sigma}$ from section \ref{Dets}, where
\begin{equation}
\label{sgnch2}
J = \{1,3,4\}, \quad \sigma(d_1,d_2,d_3,d_4,d_5) =  (-1)^{d_3 d_2 + d_5(d_4-d_5)} .
\end{equation}
The reason for this choice of the sign is that the above exact sequence can be seen as the exact sequence associated to the composition of two Fredholm operators. Indeed, let
\[
R = R_{X,(V,W)} : V \oplus W \hookrightarrow (X+V) \oplus W
\]
be the inclusion, and let
\[
T = T_{X,(V,W)} : (X+V) \oplus W \rightarrow H, \quad (v,w) \mapsto v-w,
\]
be the difference mapping. Then $R$ and $T$ are Fredholm operators, and the exact sequence (\ref{exa}) associated to their composition is precisely (\ref{es1}):
\begin{eqnarray*}
\xymatrix@C=14pt{0 \ar[r] & 0 \ar[r] \ar@{=}[d] & V\cap W \ar[r] \ar@{=}[d] & (X+V) \cap W \ar[r] \ar@{=}[d] & \frac{X+V}{V} \cong X \ar[r] \ar@{=}[d] & \frac{H}{V+W} \ar[r] \ar@{=}[d] & \frac{H}{X+V+W} \ar[r] \ar@{=}[d] & 0 \\ 0 \ar[r] & \ker R \ar[r] & \ker TR \ar[r]& \ker T \ar[r] & \coker R \ar[r] & \coker TR \ar[r] & \coker T \ar[r] & 0.}
\end{eqnarray*} 
The choice of the sign in (\ref{sgnch2}) is the same as the one in (\ref{sgnch1}).

\paragraph{Analyticity of the sum.} Consider
the set
\[  
\mathcal{S} (H) := \set{(X,(V,W))\in \Bigl(\bigcup_{n\in
    \N}\mathrm{Gr}_{n,\infty}(H) \Bigr) \times \mathrm{Fp}(H)}{X \cap V=(0)}, \\
\]
and the analytic map
\[
s: \mathcal{S}(H) \rightarrow \mathrm{Fp}(H), \quad (X,(V,W)) \mapsto (X+V,W).
\]
The space $\mathcal{S}(H)$ is the base space of the analytic line bundle
$\Det(\mathcal{S}) \rightarrow \mathcal{S}(H)$, whose fiber at
$(X,(V,W))$ is $\Det(X) \otimes \Det(V,W)$.
The collection of the isomorphisms $S(X,(V,W))$ defines a bundle morphism
\[
S:\Det(\mathcal{S}(H)) \rightarrow \Det(\mathrm{Fp}(H)),
\]
which lifts the map $s$. We can use the fact that the composition of Fredholm operators has an analytic lift to the determinant bundles to show that the above bundle morphism is analytic. Indeed, the diagram
\[
\xymatrix@C=60pt{ \Det(X) \otimes \Det(V,W) \ar@{=}[d] \ar[r]^{S(X,(V,W))} & \Det(X+V,W) \ar@{=}[d] \\ \Det(R)^* \otimes \Det (TR) \ar[r] & \Det(T)}
\]
commutes, and the lower isomorphism is induced by the inverse of the composition lift tensorized by the identity on $\Det(R)^*$.

\paragraph{Associativity of the sum.} The operation $S$ is
associative, meaning that if $X,Y\in \mathrm{Gr}(H)$ are finite dimensional,
$(V,W)\in \mathrm{Fp}(H)$, and $X\cap Y = (X+Y)\cap V = (0)$, then the diagram
\begin{equation*}
\begin{CD}
\Det(X) \otimes \Det(Y) \otimes \Det(V,W) @>{\mathrm{id}\otimes
  S(Y,(V,W))}>> \Det(X)\otimes \Det(Y+V,W) \\ @V{S(X,Y) \otimes
  \mathrm{id}}VV @VV{S(X,(Y+V,W))}V \\ \Det(X+Y) \otimes \Det(V,W)
  @>{S(X+Y,(V,W))}>> \Det(X+Y+V,W) 
\end{CD}
\end{equation*} 
commutes. This follows from the associativity property for the composition lift of Fredholm operator and from the commutativity of the diagram
\[
\xymatrix@C=100pt{V\oplus W \ar^{R_{Y,(V,W)}}[r] \ar_{R_{X+Y,(V,W)}}[d] & (Y + V) \oplus W \ar[ld]_{R_{X,(Y+V,W)}} \ar^{T_{Y,(V,W)}}[d] \\
(X+Y+V)\oplus W \ar[r]_{{T_{X+Y,(V,W)} = T_{X,(Y+V,W)}}} & H.}
\]
 
\section{Orientation bundles} 

The determinant bundle over the space of real Fredholm operators is often used in global analysis to orient a finite dimensional manifold, which is obtained as the zero set - or more generally as the inverse image of some other finite dimensional manifold - of some nonlinear real Fredholm map. See for instance \cite{fh93}.  Similarly, the determinant bundle over the space of Fredholm pairs can be used to orient finite dimensional submanifolds which are obtained as transverse intersections of infinite dimensional manifolds. See \cite{ama05} for applications of the objects introduced in this section to infinite dimensional Morse theory.

The advantage of using determinant bundles for these kind of problems lies in the associativity property. See also \cite{pej07} for other approaches to the orientation question. 

\paragraph{The orientation bundle over the space of Fredholm pairs.}
Throughout this section we assume that the Hilbert space $H$ is real, so that the determinant bundle over $\mathrm{Fp}(H)$ is a real line bundle. Its $\Z_2$-reduction 
defines a double covering 
\[
\mathrm{Or}(\mathrm{Fp}(H)) \rightarrow \mathrm{Fp}(H),
\]
whose fiber at $(V,W)\in \mathrm{Fp}(H)$ is the quotient
\[
\mathrm{Or}(V,W):= \Det(V,W) \setminus \{0\}/\sim, \quad \mbox{where } \xi\sim \eta
\mbox{ iff } \xi=\lambda \eta \mbox{ with } \lambda>0.
\]
Such a $\Z_2$-bundle will be said the {\em orientation bundle over}
$\mathrm{Fp}(H)$. Since the determinant bundle is nontrivial over each connected
component of $\mathrm{Fp}(H)$ - except the two trivial ones - 
so is the orientation bundle. 
Since in the real case the fundamental group of each
connected component of $\mathrm{Fp}^*(H)$ is $\Z_2$ (see Theorem \ref{hhoomm}), 
$\mathrm{Or}(\mathrm{Fp}^*(H)) \rightarrow \mathrm{Fp}^*(H)$ is the universal covering of $\mathrm{Fp}^*(H)$. 
If $(V,W)\in \mathrm{Fp}(H)$, $X\in \mathrm{Gr}(H)$ is finite dimensional, and $X\cap V=
(0)$, then any orientation of two of
\[
X, \quad (V,W), \quad (X + V,W),
\]
determines, by the isomorphism $S(X,(V,W))$, 
an orientation of the third one. By the
properties of $S$, this way of
summing orientations is continuous and associative. 

\paragraph{Co-orientations.} 
Let $V,W\in \mathrm{Gr}(H)$, $W$ a compact perturbation of $V$. If
$V'$ is a linear complement of $V$ in $H$, then $(W,V')$ is a Fredholm
pair (by Proposition  \ref{ugo}). Moreover, the set
\[
\mathcal{C}(W,V) := \set{(W,V')}{V' \mbox{ is a linear complement of } V
  \mbox{ in } H}
\]
is contractible, being homeomorphic to the Banach space 
$\mathrm{L}(V^{\perp},V)$. In particular, the restriction of the orientation
bundle of Fredholm pairs to $\mathcal{C}(W,V)$ is trivial, and we can give
the following:

\begin{defn}
Let $W$ be a compact perturbation of $V$.
A {\em co-orientation} of $(W,V)$ is the choice of one of the two 
continuous sections of the trivial double covering
\[
\mathrm{Or}(\mathrm{Fp}(H))|_{\mathcal{C}(W,V)}.
\]
\end{defn} 

If $H$ is endowed with a preferred inner product, we can identify a
co-orientation of $(W,V)$ with an orientation of
$(W,V^{\perp})$, hence with an
orientation of the finite dimensional space $(W\cap V^{\perp}) \oplus
(W^{\perp} \cap V)$. It follows that the concept of co-orientation is symmetric, meaning that
a co-orientation of $(W,V)$ canonically induces a co-orientation of
$(V,W)$.  

The set of the two co-orientations of $(W,V)$ is denoted by
$\mathrm{co\-Or}(V,W)$. These sets are the fibers of the 
co-orientation bundle $\mathrm{co\-Or}(H)$, a double covering of the space
\[
\set{(W,V)\in \mathrm{Gr}(H) \times \mathrm{Gr}(H)}{W \mbox{ is a
    compact perturbation of } V}.
\]
If we restrict this double covering to the Grassmannian $\mathrm{Gr_c}(V,H)$ of compact perturbations of a fixed $V\in \mathrm{Gr}_{\infty,\infty}(H)$, we obtain a nontrivial $\Z_2$-bundle
\[
\mathrm{co\-Or}(\mathrm{Gr_c}(V,H))\rightarrow \mathrm{Gr_c}(V,H),
\]
called the {\em co-orientation bundle over} $\mathrm{Gr_c}(V,H)$.
By Theorem \ref{hgpc}, this is the universal covering of
$\mathrm{Gr_c}(V,H)$. 

\paragraph{Induced orientations and co-orientations.}
Now let $(V,Z)$ be a Fredholm pair, and let $W$ be a compact
perturbation of $V$. We know from Proposition \ref{ugo} that $(W,Z)$ is also a Fredholm pair. Let us show how the choice of 
two among the following three objects
\begin{equation}
\label{three}
o_{(V,Z)} \in \mathrm{Or}(V,Z), \quad 
o_{(W,Z)} \in \mathrm{Or}(W,Z), \quad
co_{(W,V)} \in \mathrm{co\-Or}(W,V),
\end{equation}
determines the third one. The argument mimics the proof of Proposition \ref{ugo}.
  
The pair $(W,Z)$ is Fredholm if and only if the operator 
$P_{Z^{\perp}}|_W \in \mathrm{L}(W,Z^{\perp})$ is Fredholm, and 
\[
\Det(P_{Z^{\perp}}|_W) = \Det(W\cap Z) \otimes
\Det(Z^{\perp}/P_{Z^{\perp}} W) \cong  \Det(W\cap Z) \otimes
\Det(H/(W+Z)) = \Det (W,Z),
\]
where we have used the fact that $P_{Z^{\perp}} W + Z = W + Z$. Since $W$ is a compact perturbation of $V$, we can apply the above facts to the Fredholm pair $(V,W^{\perp})$, obtaining that $P_W|_V\in \mathrm{L}(V,W)$ is Fredholm and 
\[
\Det(P_W|_V) = \Det(V,W^{\perp}).
\]
By (\ref{idecomp}) we have
\[
P_{Z^{\perp}}|_V = P_{Z^{\perp}}|_W \circ P_W|_V + P_{Z^{\perp}} \circ (P_V - P_W)|_V.
\]
The last term is a compact operator, so $P_{Z^{\perp}}|_V$ is a compact perturbation of
the composition $T:=P_{Z^{\perp}}|_W \circ P_W|_V$. Using the composition lift isomorphism given by (\ref{compo}), we obtain an isomorphism
\[
\Det (T) \cong \Det (P_W|_V) \otimes \Det (P_{Z^{\perp}}|_W) \cong
\Det (V,W^{\perp}) \otimes \Det (W,Z).
\]
Since $\mathrm{L_c}(V,Z^{\perp})$ is simply connected (it is actually
contractible), an orientation of $\Det (T)$ determines an orientation
of $\Det (T')$ for each compact perturbation $T'$ of $T$. In
particular, an orientation of $\Det (T)$ determines an orientation of
\[
\Det (P_{Z^{\perp}}|_V) \cong \Det (V,Z).
\]
Therefore an orientation of
\[
\Det (V,W^{\perp}) \otimes \Det (W,Z)
\] 
determines an orientation of $\Det (V,Z)$. We conclude that the choice of 
two among the three objects in (\ref{three}) determines the third one. 

If we exchange the role of $V$ and $W$ in the above construction, we still get the same way of inducing orientations. This follows from the fact that the diagram
\[
\xymatrix{V \ar@/_/_{P_W|_V}[rr] \ar[dr]_{P_{Z^{\perp}}|_V} & & W \ar@/_/_{P_V|_W}[ll] \ar[ld]^{P_{Z^{\perp}}|_W} \\ & Z^{\perp} & }
\]
commutes up to compact perturbations. Moreover, the way of inducing orientations does not depend on the choice of the Hilbert product.

Here is a typical application of this construction: We are given a
Hilbert manifold $\mathscr{M}$ with a preferred linear subbundle
$\mathscr{V}$ of the tangent bundle $T\mathscr{M}$. Then we have two
submanifolds $\mathscr{W}$ and $\mathscr{Z}$ of $\mathscr{M}$, such
that for every $p\in \mathscr{W}$ the tangent space of $\mathscr{W}$
at $p$ is a compact perturbation of $\mathscr{V}_p$, whereas for every
$p\in \mathscr{Z}$ the pair $(T_p\mathscr{Z},\mathscr{V}_p)$ is
Fredholm. Then a co-orientation of $T\mathscr{Z}$ with respect to
$\mathscr{V}$ and an orientation of $(T\mathscr{Z},\mathscr{V})$
determine an orientation of the Fredholm pair $(T\mathscr{W},
T\mathscr{Z})$, at every point in the intersection $\mathscr{W}\cap
\mathscr{Z}$. If moreover this intersection is transverse, such an
orientation is an orientation of the finite dimensional manifold
$\mathscr{W}\cap \mathscr{Z}$. See \cite{ama05} for an application of
these concepts to infinite dimensional Morse theory.         

\paragraph{Associativity.}  Let $V,W,Z$ be as before, and consider a compact perturbation $Y$ of $Z$. Then also $(V,Y)$ and $(W,Y)$ are Fredholm pairs, and we can consider six elements
\begin{eqnarray*}
& a\in \mathrm{Or}(V,Z), \quad b\in \mathrm{Or}(V,Y), \quad c \in \mathrm{Or}(W,Z), \quad d\in \mathrm{Or}(W,Y), & \\ & e\in \mathrm{co\-Or}(V,W), \quad f\in \mathrm{co\-Or}(Y,Z),&
\end{eqnarray*}
to be associated with the edges of the tetrahedron with faces $\{a,c,e\}$, $\{a,b,f\}$, $\{b,d,e\}$, $\{c,d,f\}$:
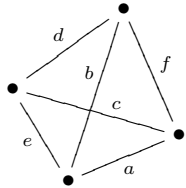
\[
\xymatrix@R=7pt@C=10pt{& & \bullet \ar@{-}[lldd]_d \ar@{-}[ldddd]_(.4){b} \ar@{-}[rddd]^f & \\ & & & \\ \bullet \ar@{-} [rdd]_e \ar@{-} '[rrrd]^(.6){c}  & & & \\ & & & \bullet \ar@{-}[lld]^{a} \\ & \bullet & & }
\]  
Then associativity can be stated in this way: If we are given five of these elements in such a way that the two triplets corresponding to the complete faces in the above tetrahedron are compatible, then there exists a unique element to be put on the remaining edge which makes also the triplets corresponding to the remaining faces compatible.  

In order to prove this fact, we may choose the Hilbert product in such a way that $Y$ - hence also $Z$ - is a compact perturbation of the orthogonal complement of $V$ - hence also of $W$. The choice of the orientations for the triplets associated to the four faces is determined by the following four compositions of Fredholm operators
\begin{eqnarray*}
V \stackrel{P_W}{\longrightarrow} W  \stackrel{P_{Y^{\perp}}}{\longrightarrow} Y^{\perp}, \quad V \stackrel{P_W}{\longrightarrow} W  \stackrel{P_{Z^{\perp}}}{\longrightarrow} Z^{\perp}, \\
Z \stackrel{P_Y}{\longrightarrow} Y  \stackrel{P_{V^{\perp}}}{\longrightarrow} V^{\perp}, \quad Z \stackrel{P_Y}{\longrightarrow} Y  \stackrel{P_{W^{\perp}}}{\longrightarrow} W^{\perp}.
\end{eqnarray*}
By the canonical isomorphism between the determinant line of  a Fredholm pair and the determinant line of the pair consisting of the orthogonal complements, we may replace the last two compositions by the compositions
\[
Z^{\perp} \stackrel{P_{Y^{\perp}}}{\longrightarrow} Y^{\perp}  \stackrel{P_V}{\longrightarrow} V, \quad Z^{\perp} \stackrel{P_{Y^{\perp}}}{\longrightarrow} Y^{\perp}  \stackrel{P_W}{\longrightarrow} W.
\]
Then the claim follows from the fact that the diagram
\[
\xymatrix@C40pt{V \ar[r]^{P_W} & W \ar[d]^{P_{Z^{\perp}}} \ar@/_/[dl]_{P_{Y^{\perp}}} \\ Y^{\perp} \ar[u]^{P_V} \ar@/_/[ru]_{P_W}& Z^{\perp} \ar[l]^{P_{Y^{\perp}}} }
\]
commutes up to compact perturbations, together with the associativity property for the composition of Fredholm operators lifted to the determinant bundles.

\paragraph{Final remarks.} Let us consider again the case of a Fredholm pair $(V,Z)$ and of a compact perturbation $W$ of $V$, and let us make some comments on the argument following (\ref{three}).

Notice that the compact operator $P_{Z^{\perp}} \circ (P_V -
P_W)|_V$ vanishes if and only if $P_{W^{\perp}} V\subset Z$. In this
situation, $P_{Z^{\perp}}|_V$ coincides with the composition $T$, and
there is a natural isomorphism already at the level of determinants,
\[
\Det(V,Z) \cong \Det (V,W^{\perp}) \otimes \Det (W,Z).
\]
This isomorphism is induced by the exact sequence (\ref{exa}), which
in this case is
\[
0 \rightarrow V\cap W^{\perp} \rightarrow V\cap Z
\stackrel{P_W}{\rightarrow} W\cap Z \rightarrow \frac{W}{P_W V}
\stackrel{P_{Z^{\perp}}}{\rightarrow} \frac{Z^{\perp}}{P_{Z^{\perp}} P_W
  V} \rightarrow \frac{Z^{\perp}}{P_{Z^{\perp}W}}  \rightarrow 0.
\]
In particular, if $X$ is a finite dimensional subspace such that
$X\cap V=(0)$ and $W=X+V$, the above exact sequence reduces to
\[
0 \rightarrow 0 \rightarrow V\cap Z \rightarrow (X+V)\cap Z
\rightarrow \frac{X+V}{V} \cong X \rightarrow
\frac{Z^{\perp}}{P_{Z^{\perp}} V} \rightarrow
\frac{Z^{\perp}}{P_{Z^{\perp}} (X+V)} \rightarrow 0.
\]
Since $P_{Z^{\perp}} V + Z = V+Z$, we have natural isomorphisms
\[
\frac{Z^{\perp}}{P_{Z^{\perp}} V} \cong \frac{H}{V+Z}, \quad 
\frac{Z^{\perp}}{P_{Z^{\perp}} (X + V)} \cong \frac{H}{X+V+Z}.
\]
Therefore, the above exact sequence is the one inducing the
isomorphism
\[
S(X,(V,W)) : \Det(X) \otimes \Det(V,Z) \cong \Det(X+V,Z).
\]
Hence, the way of inducing orientations presented in this section agrees with the sum construction of section \ref{odfp}. 

\renewcommand{\thesection}{\Alph{section}}
\setcounter{section}{0}

\section{Appendix - Sections}

A continuous linear surjective map $T:E\rightarrow F$ between Banach
spaces has a continuous linear global section, that is
a continuous linear map $S: F \rightarrow E$ such that
$TS=I_F$, if and only if the kernel of $T$ has a direct summand in
$E$. For instance, the quotient projection $\ell^{\infty}
\rightarrow \ell^{\infty}/c_0$, where $\ell^{\infty}$ denotes the
Banach space of bounded sequences, and $c_0$ denotes the closed
subspace of infinitesimal sequences, has no continuous linear right
inverse. Indeed, $c_0$ does not have a direct summand in
$\ell^{\infty}$ (see footnote \ref{fnote}).  

However, Bartle and Graves \cite{bg52} have shown that
it is always possible to find a nonlinear continuous right
inverse of $T$. In general, such a global section will not be
Gateaux differentiable at any point (otherwise its differential 
would be a linear
global section). The existence of a continuous global section could be
seen as a consequence of more general selection theorems by Michael (see
\cite{mic59} and references therein, or \cite{pal66}, Theorems 10,11). 
The aim of this appendix is to show a direct proof of this fact.

\begin{prop}
\label{propa}
Let $T:E\rightarrow F$ be a continuous linear surjective map between
Banach spaces. Then there exists a continuous map 
$s:F\mapsto E$ which is a global section of $T$.
\end{prop}

\proof
Since $T$ is onto, by the open mapping theorem the quotient norm
\[
\|x\|:= \inf\set{\|y\|_E}{Ty=x}
\]
is equivalent to the norm of $F$. We will use such a quotient norm to
define the uniform norm $\|\cdot\|_{\infty}$ of a $F$-valued map. 
Set $S:= \set{x\in F}{\|x\|=1}$.

\medskip

\noindent {\sc Claim.} If $\Psi:S\rightarrow F$ is continuous and bounded,
there exists $\Phi:S \rightarrow E$ continuous such that
$\|\Phi\|_{\infty}\leq \|\Psi\|_{\infty}$ and $\|T\Phi -
\Psi\|_{\infty}\leq (1/2) \|\Psi\|_{\infty}$. 

\medskip

For every $x\in S$, let $U(x)$ be an open neighborhood of $x$ such
that $\|\Psi(x^{\prime}) - \Psi(x)\|\leq \|\Psi\|_{\infty}/4$ for
every $x^{\prime}\in U(x)$. Let $\{V_i\}_{i\in I}$ be a locally finite
open refinement of $\{U(x)\}_{x\in S}$, and let $\{\varphi_i\}_{i\in
I}$ be a partition of unity subordinated to it.  
For every $i\in I$, let $x_i\in V_i$, and choose $y_i\in E$ such that
$Ty_i = (2/3) \Psi(x_i)$, $\|y_i\|_E \leq \|\Psi\|_{\infty}$. Define
$\Phi:S \rightarrow E$ as
\[
\Phi(x) = \sum_{i\in I} \varphi_i(x) y_i.
\]
Then $\|\Phi\|_{\infty}\leq \sup_{i\in I} \|y_i\|_E \leq
\|\Psi\|_{\infty}$, and if $x\in S$,
\[
T\Phi(x)-\Psi(x) = \sum_{i\in I} \varphi_i(x) [Ty_i - \Psi(x)] =
\frac{2}{3} \sum_{i\in I} \varphi_i(x)
[\Psi(x_i) - \Psi(x)] - \frac{1}{3} \Psi(x),
\]
hence
\[
\|T\Phi(x) - \Psi(x)\| \leq \frac{2}{3} \sup_{\substack{i\in I \\
x\in V_i}} \|\Psi(x_i) - \Psi(x)\| + \frac{1}{3}
    \|\Psi\|_{\infty} \leq \frac{1}{6} \|\Psi\|_{\infty} + \frac{1}{3}
    \|\Psi\|_{\infty} = \frac{1}{2} \|\Psi\|_{\infty},
\]
proving the claim.

\medskip

Applying the above claim iteratively, starting from the identity map,
we obtain a sequence of continuous maps $\Phi_n:S\rightarrow E$,
$n\geq 0$, such that 
\[
\|\Phi_n\|_{\infty} \leq 2^{-n} , \quad \Bigl\|I - T \sum_{k=0}^{n-1}
\Phi_k\Bigr\|_{\infty} \leq 2^{-n}.
\]
Therefore, the series $\sum_{n=0}^{\infty} \Phi_n$ converges uniformly
on $S$ to a map which lifts the identity, and the map
\[
s(x) := \|x\| \sum_{n=0}^{\infty} \Phi_n\left( \frac{x}{\|x\|}
\right), \quad s(0)=0,
\]
is a global section of $T$.
\qed   

\begin{rem}
From the existence of the global section $s$, it follows easily that
every continuous linear surjective map $T:E\rightarrow F$ admits a
natural structure of a trivial $C^0$ vector bundle over $F$, and the
map $x\mapsto \|x-s(Tx)\|_E$ is a Finsler structure on such a vector
bundle. 
\end{rem}

\section{Appendix - Associativity}

Consider the Fredholm operators $T_i \in \mathrm{Fr}(H_i,H_{i+1})$, $1\leq i \leq 3$. We wish to prove that the diagram 
\begin{equation}
\label{diago2}
\begin{CD}
\Det(T_1) \otimes \Det(T_2) \otimes \Det(T_3) @>{\mathrm{id}\otimes
  \phi_{T_2,T_3}}>> \Det(T_1) \otimes \Det(T_3 T_2) \\ @VV{\phi_{T_1,T_2} \otimes
  \mathrm{id}}V @VV{\phi_{T_1,T_3 T_2}}V \\ \Det(T_2 T_1) \otimes \Det(T_3)
  @>{\phi_{T_2 T_1,T_3}}>> \Det(T_3 T_2 T_1)
\end{CD}
\end{equation} 
commutes. 

Given $1\leq i < j \leq 4$, set $T_{ij} = T_{j-1} 
\circ \dots \circ T_i \in \mathrm{Fr}(H_i,H_j)$. Using the left and right
action of the general linear group,
we are reduced to consider the following situation:
\[
H_i = \bigoplus_{h\leq i \leq k} H_{hk}, \quad T_{ij} =
\bigoplus_{\substack{h\leq i \\ k \geq j}} I_{H_{hk}},
\]
where the space $H_{hk}$, $1\leq h\leq k\leq 4$, has finite
dimension $d_{hk}$, except for $H_{14}$ which is infinite
dimensional. Then
\[
\ker T_{ij} = \bigoplus_{\substack{h\in [1,i]\\ k\in [i,j-1]}} H_{hk},
\quad
\coker T_{ij} = \bigoplus_{\substack{h\in [i+1,j]\\ k\in [j,4]}} H_{hk}.
\]
Fix $1\leq i<j<\ell \leq 4$.
The exact sequence associated to the composition $T_{i\ell} =
T_{j\ell} \circ T_{ij}$ is 
\[
\xymatrix@C=7pt@R=10pt{0 \ar[r] & \ker T_{ij} \ar@{=}[d] \ar[r] & \ker T_{i\ell}
  \ar@{=}[d] \ar[r] & \ker T_{j\ell} \ar@{=}[d] \ar[r] & \coker T_{ij}  
\ar@{=}[d] \ar[r] & \coker T_{i\ell} \ar@{=}[d] \ar[r] & \coker
T_{j\ell} \ar@{=}[d] \ar[r] & 0 \\
0 \ar[r] & \displaystyle{\bigoplus_{\substack{h\in [1,i]\\ k\in [i,j-1]}} H_{hk}}
\ar[r] & \displaystyle{\bigoplus_{\substack{h\in [1,i]\\ k\in [i,\ell-1]}} H_{hk}}
\ar[r] & \displaystyle{\bigoplus_{\substack{h\in [1,j]\\ k\in [j,\ell-1]}} H_{hk}}
\ar[r] & \displaystyle{\bigoplus_{\substack{h\in [i+1,j]\\ k\in [j,4]}} H_{hk}}
\ar[r] & \displaystyle{\bigoplus_{\substack{h\in [i+1,\ell]\\ k\in
      [\ell,4]}} H_{hk} }
\ar[r] & \displaystyle{\bigoplus_{\substack{h\in [j+1,\ell]\\ k\in
      [\ell,4]}} H_{hk} }
\ar[r] & 0.}
\]
Denote by $\theta_{hk}$ a generator of $\Det(H_{hk})$, for $(h,k) \neq
(1,4)$, and given two sets of consecutive integers $I,J\subset
\{1,2,3,4\}$, set 
\[
\Theta_I^J = \bigwedge_{h\in I} \bigwedge_{k\in J} \theta_{hk},
\]
where we are considering the standard order in each wedge product.
Let us show that the composition morphism
\[
\phi_{T_{ij},T_{j\ell}} : \Det(T_{ij}) \otimes \Det(T_{j\ell})
\rightarrow \Det(T_{i\ell})
\]
equals
\begin{equation}
\label{sig}
\phi_{T_{ij},T_{j\ell}}: {\Theta_{[1,i]}^{[i,j-1]}} \otimes
{\Theta_{[i+1,j]}^{[j,4]}}^* \otimes {\Theta_{[1,j]}^{[j,\ell-1]}} \otimes
{\Theta_{[j+1,\ell]}^{[\ell,4]}}^* \mapsto 
(-1)^{\sigma(i,j,\ell)} \, {\Theta_{[1,i]}^{[i,\ell-1]}} \otimes
{\Theta_{[i+1,\ell]}^{[\ell,4]}}^*,
\end{equation}
where the coefficient $\sigma(i,j,\ell)$ is to be determined.   
Notice that:

\begin{enumerate}

\item If we want to change the generator
  $\Theta_{[1,i]}^{[i,\ell-1]}$ of $\Det(\ker
  T_{i\ell})$ into $\Theta_{[1,i]}^{[i,j-1]} \wedge 
  \Theta_{[1,i]}^{[j,\ell-1]}$, we have to exchange the position of
  $\theta_{hk}$ and $\theta_{h'k'}$ if and only if $h<h'$, $k\in
  [j,\ell-1]$, and $k'\in [i,j-1]$. Therefore,
\[
\Theta_{[1,i]}^{[i,\ell-1]} = (-1)^{\sigma_1(i,j,\ell)}\,
\Theta_{[1,i]}^{[i,j-1]} \wedge \Theta_{[1,i]}^{[j,\ell-1]},
\]
where
\[
\sigma_1(i,j,\ell) = \sum_{\substack{1 \leq h < h' \leq i\\ k \in
    [j,\ell-1] \\ k' \in [i,j-1]}} d_{hk} d_{h'k'}.
\]

\item The generator $\Theta_{[1,j]}^{[j,\ell-1]}$ of $\Det(\ker
T_{j\ell})$ can be rewritten as $\Theta_{[1,i]}^{[j,\ell-1]} \wedge
\Theta_{[i+1,j]}^{[j,\ell-1]}$. Therefore, 
\[
\Theta_{[1,j]}^{[j,\ell-1]} = (-1)^{\sigma_2(i,j,\ell)}\, 
\Theta_{[i+1,j]}^{[j,\ell-1]} \wedge \Theta_{[1,i]}^{[j,\ell-1]},
\]
where
\[
\sigma_2(i,j,\ell) = \sum_{\substack{h \in [1,i]\\
h'\in [i+1,j] \\ k,k' \in [j,\ell-1]}} d_{hk} d_{h'k'}.
\]

\item Arguing as in (i), the generator $\Theta_{[i+1,j]}^{[j,4]}$ of
  $\Det(\coker T_{ij})$ can be rewritten as
\[
\Theta_{[i+1,j]}^{[j,4]} = (-1)^{\sigma_3(i,j,\ell)}\,
\Theta_{[i+1,j]}^{[j,\ell-1]} \wedge \Theta_{[i+1,j]}^{[\ell,4]},
\]
where
\[
\sigma_3(i,j,\ell) = \sum_{\substack{i+1 \leq h < h' \leq j\\ k \in
    [\ell,4] \\ k' \in [j,\ell-1]}} d_{hk} d_{h'k'}.
\]

\item As in (ii), the generator $\Theta_{[i+1,\ell]}^{[\ell,4]}$ of $\Det(\coker
T_{i\ell})$ can be rewritten as
\[
\Theta_{[i+1,\ell]}^{[\ell,4]} = (-1)^{\sigma_4(i,j,\ell)} \, 
\Theta_{[j+1,\ell]}^{[\ell,4]} \wedge \Theta_{[i+1,j]}^{[\ell,4]},
\]
where
\[
\sigma_4(i,j,\ell) = \sum_{\substack{h \in [i+1,j]\\
\\ h'\in [j+1,\ell] \\
k,k' \in [\ell,4]}} d_{hk} d_{h'k'}.
\]

\end{enumerate}

Then formula (\ref{compo}) implies formula (\ref{sig}) with
\[
\sigma(i,j,\ell) = \sigma_0(i,j,\ell) + \sigma_1(i,j,\ell) +
\sigma_2(i,j,\ell) + 
\sigma_3(i,j,\ell) + \sigma_4(i,j,\ell),
\]
where
\[
\sigma_0(i,j,\ell) := (\dim \ker T_{ij} + \dim \ker T_{i\ell})(\dim
\coker T_{j\ell} + \dim \coker T_{i\ell}) = \sum_{\substack{h\in
    [1,i]\\ k \in [j,\ell-1]\\ h' \in [i+1,j]\\ k' \in [\ell,4]}}
d_{hk} d_{h'k'} \quad (\mbox{mod } 2).
\]
Formula (\ref{sig}) implies that the diagram (\ref{diago2}) commutes up
to the sign
\[
(-1)^{\sigma(1,2,3)+\sigma(1,2,4) + \sigma(1,3,4) + \sigma(2,3,4)},
\]
so we have to show that the above exponent vanishes modulo 2. The
following computations are modulo 2. We have
\[
\sigma_0(i,j,\ell) + \sigma_2(i,j,\ell) = \sum_{\substack{h\in [1,i]\\
    k\in [j,\ell-1],\\ h'\in [i+1,j], \\ k'\in [j,4]}} d_{hk}
d_{h'k'},
\]
and summing over the four possibilities for $(i,j,\ell)$ we obtain
\begin{equation}
\label{ass1}
(\sigma_0+\sigma_2)(1,2,3) + (\sigma_0+\sigma_2)(1,2,4) +
(\sigma_0+\sigma_2)(1,3,4) + (\sigma_0+\sigma_2)(2,3,4) = d_{13}d_{22}
+ d_{23} d_{33} + d_{23} d_{34}.
\end{equation}
The range of the sum defining $\sigma_1(i,j,\ell)$ is always empty, except for
the case $i=2$, $j=3$, $\ell=4$, so we have
\begin{equation}
\label{ass2}
\sigma_1(1,2,3) + \sigma_1(1,2,4) +
\sigma_1(1,3,4) + \sigma_1(2,3,4) =
\sigma_1(2,3,4) = d_{13} d_{22}.
\end{equation}
The range of the sum defining $\sigma_3(i,j,\ell)$  is always empty, except for
the case $i=1$, $j=3$, $\ell=4$, so we have
\begin{equation}
\label{ass3}
\sigma_3(1,2,3) + \sigma_3(1,2,4) +
\sigma_3(1,3,4) + \sigma_3(2,3,4) =
\sigma_3(1,3,4) = d_{24} d_{33}.
\end{equation}   
Finally, taking the sum of $\sigma_4(i,j,\ell)$ over 
 the four possibilities for $(i,j,\ell)$ we obtain
\begin{equation}
\label{ass4}
\sigma_4(1,2,3) + \sigma_4(1,2,4) +
\sigma_4(1,3,4) + \sigma_4(2,3,4) = d_{23}d_{33}
+ d_{23} d_{34} + d_{24} d_{33}.
\end{equation}
Taking the sum of (\ref{ass1}), (\ref{ass2}), (\ref{ass3}), and
(\ref{ass4}), we find
\[
\sigma(1,2,3)+\sigma(1,2,4) + \sigma(1,3,4) + \sigma(2,3,4) = 0,
\]
as wished.

\providecommand{\bysame}{\leavevmode\hbox to3em{\hrulefill}\thinspace}
\providecommand{\MR}{\relax\ifhmode\unskip\space\fi MR }
\providecommand{\MRhref}[2]{%
  \href{http://www.ams.org/mathscinet-getitem?mr=#1}{#2}
}
\providecommand{\href}[2]{#2}

\end{document}